\documentclass[11pt]{paper}
\usepackage{amssymb}
\usepackage{amsmath}
\usepackage{amsthm}
\usepackage{array}
\usepackage{graphicx}
\usepackage{subfigure}
\usepackage{epstopdf}
\usepackage{cases}
\usepackage{color}
\usepackage{multirow}
\usepackage{picinpar}
\usepackage[colorlinks,
            linkcolor=red,
            anchorcolor=green,
            citecolor=blue
            ]{hyperref}
\usepackage{enumitem}
\usepackage{fourier}
\usepackage{bookmark}
\begin{document}

\newcommand{\defi}{\stackrel{\Delta}{=}}
\newcommand{\A}{{\cal A}}
\newcommand{\B}{{\cal B}}
\newcommand{\U}{{\cal U}}
\newcommand{\G}{{\cal G}}
\newcommand{\cZ}{{\cal Z}}
\newcommand\one{\hbox{1\kern-2.4pt l }}
\newcommand{\Item}{\refstepcounter{Ictr}\item[\left(\theIctr\right)]}
\newcommand{\QQ}{\hphantom{MMMMMMM}}

\newtheorem{Theorem}{Theorem}[section]
\newtheorem{Lemma}{Lemma}[section]
\newtheorem{Corollary}{Corollary}[section]
\newtheorem{Remark}{Remark}[section]
\newtheorem{Example}{Example}[section]
\newtheorem{Proposition}{Proposition}[section]
\newtheorem{Property}{Property}[section]
\newtheorem{Assumption}{Assumption}[section]
\newtheorem{Definition}{Definition}[section]
\newtheorem{Construction}{Construction}[section]
\newtheorem{Condition}{Condition}[section]
\newtheorem{Exa}[Theorem]{Example}
\newcounter{claim_nb}[Theorem]
\setcounter{claim_nb}{0}
\newtheorem{claim}[claim_nb]{Claim}
\newenvironment{cproof}
{\begin{proof}
 [Proof.]
 \vspace{-3.2\parsep}}
{\renewcommand{\qed}{\hfill $\Diamond$} \end{proof}}
\newcommand{\erhao}{\fontsize{21pt}{\baselineskip}\selectfont}
\newcommand{\xiaoerhao}{\fontsize{18pt}{\baselineskip}\selectfont}
\newcommand{\sanhao}{\fontsize{15.75pt}{\baselineskip}\selectfont}
\newcommand{\sihao}{\fontsize{14pt}{\baselineskip}\selectfont}
\newcommand{\xiaosihao}{\fontsize{12pt}{\baselineskip}\selectfont}
\newcommand{\wuhao}{\fontsize{10.5pt}{\baselineskip}\selectfont}
\newcommand{\xiaowuhao}{\fontsize{9pt}{\baselineskip}\selectfont}
\newcommand{\liuhao}{\fontsize{7.875pt}{\baselineskip}\selectfont}
\newcommand{\qihao}{\fontsize{5.25pt}{\baselineskip}\selectfont}
\newcounter{Ictr}
\renewcommand{\theequation}{
\arabic{equation}}
\renewcommand{\thefootnote}{\fnsymbol{footnote}}

\def\A{\mathcal{A}}

\def\C{\mathcal{C}}

\def\V{\mathcal{V}}

\def\I{\mathcal{I}}

\def\Y{\mathcal{Y}}

\def\X{\mathcal{X}}

\def\J{\mathcal{J}}

\def\Q{\mathcal{Q}}

\def\W{\mathcal{W}}

\def\S{\mathcal{S}}

\def\T{\mathcal{T}}

\def\L{\mathcal{L}}

\def\M{\mathcal{M}}

\def\N{\mathcal{N}}
\def\R{\mathbb{R}}
\def\H{\mathbb{H}}

\title{}
\author{}
\begin{center}
\topskip2cm
\LARGE{\bf Tuning parameter selection rules for nuclear norm regularized multivariate linear regression}
\end{center}

\begin{center}
\renewcommand{\thefootnote}{\fnsymbol{footnote}}Pan Shang,  Lingchen Kong, \footnote{e-mail: 18118019@bjtu.edu.cn,   konglchen@126.com, }\\
Beijing Jiaotong University, China\\
(January 18th, 2019)\\
\end{center}
\vskip4pt
\textbf{Summary:} We consider the tuning parameter selection rules for nuclear norm regularized  multivariate linear regression (NMLR) in high-dimensional setting.  High-dimensional multivariate linear regression is widely used in statistics and machine learning, and regularization technique is commonly applied to deal with the special structures in high-dimensional data. As we know, how to select the tuning parameter  is an essential issue for regularization approach and it directly affects the model estimation performance. To the best of our knowledge, there are no rules about the tuning parameter selection for NMLR from the point of view of optimization. In order to establish such rules, we study the duality theory of NMLR. Then, we claim the choice of tuning parameter for NMLR is based on the sample data and the solution of NMLR dual problem, which is  a projection on a nonempty, closed and convex set.  Moreover, based on the (firm) nonexpansiveness and the idempotence of  the projection operator, we build four tuning parameter selection rules PSR, PSRi, PSRfn and PSR+.  Furthermore,  we give a sequence of tuning parameters and the corresponding intervals for every rule, which states that the rank of the estimation coefficient matrix is no more than a fixed number for the tuning parameter in the given interval. The relationships between these rules are also discussed and PSR+ is the most efficient one to select the tuning parameter. Finally, the numerical results are reported on simulation and real data, which show that  these four tuning parameter selection rules are  valuable.
\section{Introduction}
High-dimensional multivariate linear regression is widely used in many areas, such as chemometrics, econometris, engineering, gene expression and so on. A well-known example is a breast cancer study about the influence of DNA copy number alterations on RNA transcript levels (Peng et al. 2010), which includes 172 samples, 384 DNA copy number and 654 breast cancer related RNA expressions for every sample.  Here, the predictors are 384 DNA copy number and the responses are 654 breast cancer related RNA expressions. Thus, in this instance, the prediction matrix is 172 by 384 and the response matrix is 172 by 654. In order to explore the influence of DNA copy number alterations on RNA transcript levels, a direct way is to establish the multivariate linear regression model where the coefficient matrix measure this influence. Note that the sample size is less than the number of predictors or responses, which means the data is high-dimensional. For this high-dimensional matrix data, one common assumption of the coefficient matrix  is low rank, see, e.g., Yuan et al (2007) and Negahban et al (2011).  However, the optimization problems with low rank constraint are NP-hard. The regularization technique is always used to deal with these problems, and the popular regularization is nuclear norm instead of the low-rank constraint.  Hence, our concern  is nuclear norm regularized  multivariate linear regression (NMLR) in this paper. Clearly, if response variables are univariate, NMLR degrades into the famous LASSO (Tibshirani, 1996; Chen et al., 1998).

Tuning parameter selection is an important issue for regularization approach and it affects the model estimation performance. From the perspective of prediction accuracy, tuning parameter can be chosen by cross validation and information criteria. See, e.g., Wang et al. (2007), Wang et al. (2009), Fan et al. (2013) and so on. Meanwhile, there are some screening rules for LASSO under the help of optimization techniques, which eliminate the inactive predictors by choosing the appropriate tuning parameters. For example, Fan et al. (2008) proposed the sure independence screening (SIS), which reduces dimensionality of the predictors below sample size.  The idea of SIS is  to select predictors using their correlations with the response. Ghaoui et al. (2012) constructed SAFE rules that help to eliminate predictors in LASSO, which are based on the duality theorem  in optimization. The SAFE rules never remove active predictors. Specifically, they proved that applying these tests to eliminate predictors can save time and memory in computational process.  Tibshirani et al. (2012) proposed strong rules for discarding inactive predictors under the assumption of the unit slope bound.  The strong rules screen out far more predictors than SAFE rules in practice and can be more efficient by checking Karush-Kuhn-Tucker conditions for any predictor.  Wang et al. (2015) built dual polytope projection (DPP) and the enhanced version EDPP to discard inactive predictors. They showed that EDPP had a better performance in screening out inactive predictors than SAFE rules and strong rules. These screening rules closely relate to the sparsity of the coefficient vector. Here, the sparsity means many elements of the coefficient vector are zero, which implies lots of the predictors are inactive. In the sense of the sparse solution in LASSO, screening rules are in essence the tuning parameter selection rules. By analyzing the above arguments, we know that the optimization techniques play a role in selecting the tuning parameter for LASSO. This opens a hope that we may build up tuning parameter selection rules for NMLR from the point of view of optimization.  However, the low rank of a matrix doesn't mean lots of zero elements of the matrix, but the sparsity of singular value vector. One nature question is, \emph{how to establish the tuning parameter selection rules for NMLR?}

This paper will deal with this problem and give an affirmative answer.  In order to do so, we present the dual problem of NMLR and prove that the dual solution is a projection on a nonempty, closed and convex set.  This set is not a polytope and more complex, which is different from DPP and EDPP in Wang et al. (2015). With the help of optimization technique, we show that the strong duality theorem holds on NMLR and its dual problem. This implies that the choice of the tuning parameter for NMLR is closely related to its dual solution and sample data. Secondly,  we give an estimate set for the dual solution of NMLR based on the nonexpansiveness of the projection operator. This together with the optimal conditions of the primal and dual problems, we can estimate the maximal rank of the solution of NMLR for any tuning parameter, which leads to our basic tuning parameter selection rule PSR. In the similar way, we obtain PSRi based on  the idempotence  of the projection operator, and PSRfn on the firm nonexpansiveness.  Both PSRi and PSRfn outperform  PSR, because the estimate sets for them are more accurate than that for PSR. Moreover, by combining the  idempotence and the firm nonexpansiveness of the projection operator, we continue to get the enhance version PSR+ which behaves better than PSRi and PSRfn surely. This leads to PSR+ is the best rule. Furthermore, we give a sequence of tuning parameters and the corresponding intervals for every rule, which states that the rank of the estimation coefficient matrix is no more than a fixed number for the tuning parameter  in the given interval. Thirdly, because a dual solution of NMLR is the basis of these rules, we need an efficient algorithm for solving it. Therefore, we present the detail process of the alternating direction multiplier method (ADMM)  for solving the dual problem.  Finally,  we illustrate these rules and the ADMM  on simulation and real data. The numerical results report that our tuning parameter selection rules are valuable  and PSR+ is the most efficient rule. Actually, our tuning parameter selection rules can be applicable to any efficient algorithm for solving the dual problem of NMLR.

In all, the main contributions of our paper are threefold.
\begin{enumerate}[fullwidth,itemindent=2em, topsep=0ex, parsep=0ex, label=(\roman*)]
\item We state that the tuning parameter selection rules for NMLR are connected to sample data and the dual solution, where the dual solution is a projection on a nonempty, closed and convex set.
\item We construct four tuning parameter selection rules PSR, PSRi, PSRfn and PSR+, and show their relationships based on the properties of projection operator. These rules claim sequences of tuning parameters and the corresponding intervals, where the maximal rank of the estimate coefficient matrix is given.
\item We present the detail process of  ADMM for solving the dual problem of NMLR. The numerical results on simulation and real data  demonstrate that the four tuning parameter selection rules are valuable.
\end{enumerate}

The rest of this paper is organized as follows. We present NMLR and its dual theory in Section 2. In Section 3, we show four tuning parameter selection rules based on properties of the projection operator, and give a sequence of tuning parameters and corresponding intervals for every rule. Moreover, we claim the relationships between these tuning parameter selection rules. In Section 4, we propose an efficient ADMM to solve the dual problem of NMLR and illustrate the tuning parameter selection rules are valuable in numerical study. Some conclusion remarks are given in Section 5. In appendix, we give the proof of main results.
\section{Preliminaries}
In this section, we introduce nuclear norm regularized multivariate linear regression and show its duality theory from the optimization perspective.

We begin with reviewing the statistical model of multivariate linear regression (MLR) as follows
\begin{center}
$\textbf{y}^{T}=\textbf{x}^{T}B+\epsilon^{T},$
\end{center}
where $\textbf{x}\in R^{p}$ is the prediction vector and $\textbf{y}\in R^{q}$ is the response vector, $B\in R^{p\times q}$ is the coefficient matrix which is unknown,  $\epsilon$ is a random error vector.  By sampling $n$ times, we get $$\textbf{y}_{i}^{T}=\textbf{x}_{i}^{T}B+\epsilon_{i}^{T},~~i=1,2,\cdots,n.$$
For easy of representation, let $X=(\textbf{x}_{1},\textbf{x}_{2},\cdots,\textbf{x}_{n})^{T}\in R^{n \times p}$ be the prediction matrix, $Y=(\textbf{y}_{1},\textbf{y}_{2},\cdots,\textbf{y}_{n})^{T}\in R^{n\times q}$ the response matrix and $\mathcal{E}=(\epsilon_{1},\epsilon_{2},\cdots,\epsilon_{n})^{T}\in R^{n\times q}$ the random error matrix. Then we can write MLR with n samples in matrix form
\begin{eqnarray*}
Y=XB+\mathcal{E},
\end{eqnarray*}
In this paper, we assume random error variables in $\mathcal{E}$ are all with mean 0 and standard variance $\sigma$. Usually, the least square fitting is a capable tool to estimate coefficient matrix $B$.  For high-dimensional data with n less than p or q, we assume that the coefficient matrix is low rank. In this case, regularization technique is a  popular method to deal with the special structure of the coefficient matrix. The common regularization of low-rank constraint is nuclear norm, so we focus on nuclear norm regularized multivariate linear regression (NMLR) in this paper, which is given as follows
\begin{eqnarray}
\underset{B}\min\left\{\frac{1}{2}\|Y-XB\|^{2}_{F}+\lambda||B||_{*}\right\},
\end{eqnarray}
where $\lambda\geq0$ is the tuning parameter. Clearly, when $q=1$, NMLR degrades into the famous LASSO. The  solution of NMLR (1) relies on the choice of tuning parameter $\lambda$, so we denote $B^{*}(\lambda)$ as the solution.  For any matrix $M\in R^{p\times q}$, suppose $M$ has a singular value decomposition with nondecreasing singular values  $\sigma_{1}(M)\geq \cdots \sigma_{r}(M)\geq 0$, where $r=\rm{min\left\{p, q\right\}}$ and it's used throughout this paper. There are some norms related to $M$  and these definitions are used throughout the paper. The Frobenius norm $\|\cdot\|_{F}$ is defined as $\|M\|_{F}=\sqrt{\sum_{i=1}^{p}\sum _{j=1}^{q}x_{ij}^{2}}=\sqrt{\rm tr(M^{T}M)}=\sqrt{\sigma_{1}(M)^{2}+\cdots+\sigma_{r}(M)^{2}}$. The nuclear norm $\|\cdot\|_{*}$ is  the sum of singular values, i.e., $\|M\|_{*}=\sum_{i=1}^{r}{\sigma_{i}(M)}$. The spectral norm $\|\cdot\|_{2}$ is the largest singular value, i.e., $\|M\|_{2}=\sigma_{1}(M)$.

Now, we consider about the duality theory of NMLR (1).  First, rewrite it as
\begin{align}
\underset{B,A} \min\left\{\lambda||B||_{*}+\frac{1}{2}\|A\|^{2}_{F}\right\} \quad   s.t. \quad &Y-XB-A=0.
\end{align}
Thus, we have Lagrangian function of (2)
\begin{center}
$\textit{L}\left(B,A;C\right)=\lambda||B||_{*}+\frac{1}{2}\|A\|^{2}_{F}+\left\langle C, Y-XB-A\right\rangle$.
\end{center}
where $C\in R^{n\times q}$ is a Lagrangian multiplier. We can yield the dual problem of (2)
\begin{align}
\underset{C}\min \left\{ \frac{\lambda^{2}}{2}\left\|C-\frac{Y}{\lambda}\right\|^{2}_{F}-\frac{1}{2}\|Y\|^{2}_{F}\right\} \quad
s.t. \quad &\|X^{T}C\|_{2}\leq1.
\end{align}
The detail process of duality analysis is presented in Appendix A.2. Denote the feasible area of (3) as $\Omega_{D}=\left\{C\Big| \|X^{T}C\|_{2}\leq1\right\}$. It's clear that $\Omega_{D}$ is a nonempty, closed and convex set, and the solution of (3) is
\begin{eqnarray}
C^{*}(\lambda)=P_{\Omega_{D}}\left(\frac{Y}{\lambda}\right),
\end{eqnarray}
where $P_{\Omega_{D}}(\cdot)$ denotes the projection operator on $\Omega_{D}$ (see the details in Appendix A.1). Note that $\Omega_{D}$ is not a polytope, which is different from LASSO in vector case, see, Wang et al. (2015). From the optimality conditions analysis, we have the Karush-Kuhn-Tucker (KKT) system of (2) and (3)
\begin{eqnarray}
\begin{cases}
X^{T}C \in \partial\|B\|_{*},\\
A=\lambda C,\\
Y-XB-A=0.
\end{cases}
\end{eqnarray}
If a pair $\left(B^{*}(\lambda),A^{*}(\lambda),C^{*}(\lambda)\right)$ satisfies the KKT system, it's called the KKT point of (2) and (3).  Based on the convex optimization analysis (see Appendix A.2), it holds the strong duality theorem.
\begin{Theorem}(\textbf{Strong duality theorem})\quad
Problem (2) satisfies Slater's constraint qualification and there is a KKT point $(C^{*}(\lambda),B^{*}(\lambda),A^{*}(\lambda))$ such that  the optimal values of (2) and (3) are equal, i.e.,
\begin{center}
$\lambda||B^{*}(\lambda)||_{*}+\frac{1}{2}\|A^{*}(\lambda)\|^{2}_{F}=-\left(\frac{\lambda^{2}}{2}\left\|C^{*}(\lambda)-\frac{Y}{\lambda}\right\|^{2}_{F}-\frac{1}{2}\|Y\|^{2}_{F}\right)$.
\end{center}
Here, $(B^{*}(\lambda),A^{*}(\lambda))$ is the solution of (2) and $C^{*}(\lambda)$ is the solution of (3). Moreover, $B^{*}(\lambda)$ is the solution of NMLR (1).
\end{Theorem}
\indent According to the $X^{T}C \in \partial\|B\|_{*}$ in KKT system (5) and Theorem 2.1, we easily obtain a sufficient condition for estimating the some singular values of $B^{*}(\lambda)$ being zero.
\begin{Theorem}
For any tuning parameter $\lambda$, if $\sigma_{i}(X^{T}C^{*}(\lambda))<1$, then
$\sigma_{i}(B^{*}(\lambda))=0$,
where $B^{*}(\lambda),C^{*}(\lambda)$  are solutions of (1) and (3), respectively.
\end{Theorem}
It's worth noting that the number of nonzero elements decides the sparsity of a vector and the number of nonzero singular values decides the rank of a matrix. For the purpose of selecting a tuning parameter satisfying that the rank of the solution of NMLR (1) is given,  we need to consider about its singular values, not the entries. Theorem 2.2 implies that for any fix tuning parameter, sample data and the solution of  dual problem (3) can decide whether the singular value of the solution of NMLR (1) is zero or not. Based on this, the rank of the solution of NMLR (1) can be yielded.  Note that  $\rm rank\left(X^{T}C^{*}\left(\lambda\right)\right)\leq rank (X)$, if the prediction matrix X is not full rank, neither is $X^{T}C^{*}\left(\lambda\right)$. According to Theorem 2.2, we have the following result.
\begin{Corollary}
If X is not full rank, the solution of (1) is not full rank.
\end{Corollary}
 We already know that NMLR (1) is equivalent to problem (2). Hence, we analyze the tuning parameter for NMLR (1) through the solution of (3) and Theorem 2.2. It's fortunate that the solution is a projection on a nonempty, closed and convex set,  and properties of the projection operator can help to establish the tuning parameter selection rules. 
\section{Tuning parameter selection rules}
In this section, we give four tuning parameter selection rules: PSR, PSRi, PSRfn and PSR+. These rules provide evidences to choose the tuning parameter, satisfying that the maximal rank of the solution of NMLR is decided by sample data and the solution of dual problem (3) under a fixed tuning parameter.  The differences between these rules are that they are based on different properties of the projection operator. One can see all proofs of the results in this section in Appendix B.

We start with the lower bound of the tuning parameter that enforce the solution of NMLR (1) is zero. For problem (1), we already know that $B^{*}(\lambda)=0$ if $\lambda$ is sufficiently large. The next proposition gives the lower bound of the tuning parameter $\lambda$ which guarantees $B^{*}(\lambda)=0$.
\begin{Proposition}
 $B^{*}(\lambda)=0$ is the solution of problem (1) if and only if $\lambda\geq\lambda_{max}:=\|X^{T}Y\|_{2}$.
\end{Proposition}
From Proposition 3.1, it's clear that for any $\lambda<\lambda_{max}$, $B^{*}(\lambda)\neq0$. Because we are interested in the solution $B^{*}(\lambda)$ is nonzero, $\lambda$ needs to be less than $\lambda_{max}$. Hence, we focus on the case of $\lambda$ such that $0<\lambda<\lambda_{0}$ for given $\lambda_{0}\leq\lambda_{max}$. We have the following result, which states that the rank of the solution of NMLR (1) relates to the solution of (3) and sample data.
\begin{Theorem}{\rm(\textbf{PSR})}
Assume the solution $C^{*}\left(\lambda_{0}\right)$ of (3) is known for given $\lambda_{0}$. For any $i\in \{1,\cdots,r\}$, if $\lambda<\lambda_{0}$ and
\begin{center}
$\lambda>\frac{\lambda_{0}\|X\|_{2}\|Y\|_{F}}{\lambda_{0}-\lambda_{0}\sigma_{i}\left(X^{T}C^{*}\left(\lambda_{0}\right)\right)+\|X\|_{2}\|Y\|_{F}}$,
\end{center}
then the solution of NMLR (1) satisfies  $\sigma_{j}\left(B^{*}(\lambda)\right)=0$  $(j\geq i)$. Moreover, rank$\left(B^{*}(\lambda)\right)\leq i-1$.
\end{Theorem}
Theorem 3.1 claims that if $\lambda_{0}$ is set and the solution $C^{*}(\lambda_{0})$ of (3) is easy to solve, we can  select the tuning parameter $\lambda$ such that the maximal rank of $B^{*}(\lambda)$ is certain. In general,  $C^{*}(\lambda_{0})$ may not be computed easily for a given $\lambda_{0}$. Fortunately,  for $\lambda_{0}=\lambda_{max}=\|X^{T}Y\|_{2}$, the solution $C^{*}(\lambda_{0})$ equals to $\frac{Y}{\lambda_{0}}$ from the proof of Proposition 3.1.  Now, we talk about the relationship between the rank of $B^{*}(\lambda)$ and  $X^{T}Y$. From KKT system (5) and Proposition 3.1, we obtain Proposition 3.2.
\begin{Proposition}
For any $0<\lambda<\lambda_{max}$.  If all singular values of $X^{T}Y$ are equal to a certain number, then
\begin{center}
$\rm rank$$(B^{*}(\lambda))=$$\rm rank$$(X^{T}Y).$
\end{center}
\end{Proposition}
Next, we discuss the case that $X^{T}Y$ has at least  two different singular values.
\begin{Theorem}
 Suppose $X^{T}Y$ has at least two different singular values. For $i\in \{1,\cdots,r\}$, let's define $\lambda_{i}$  as
\begin{center}
$\lambda_{i}=\frac{\|X^{T}Y\|_{2}\|X\|_{2}\|Y\|_{F}}{\|X^{T}Y\|_{2}-\sigma_{i}\left(X^{T}Y\right)+\|X\|_{2}\|Y\|_{F}}$.
\end{center}
Then, for any  $\lambda \in (\lambda_{i},\lambda_{i-1}]$ ($i\geq 2$),  the solution of NMLR (1) satisfies
\begin{center}
rank$\left(B^{*}(\lambda)\right)\leq i-1$.
\end{center}
\end{Theorem}
From the proof of Theorem 3.1, we know $C^{*}\left(\lambda\right) \in \Omega=\left\{C\Big|\|C-C^{*}\left(\lambda_{0}\right)\|_{F}\leq \left(\frac{1}{\lambda}-\frac{1}{\lambda_{0}}\right)\|Y\|_{F}\right\}$. Denote $\rho=\left(\frac{1}{\lambda}-\frac{1}{\lambda_{0}}\right)\|Y\|_{F}$ as the radius of this set, it's clear that the set of $C^{*}(\lambda)$ will be more accurate with $\rho$ decreasing. Therefore, the aim of next parts is to reach a smaller $\rho$ which directly results in improvement consequences of PSR. The tuning parameter selection rule PSR is based on the basic property of the projection operator. A nature idea is to improve the results by using the other properties of the projection operator. Before doing so,  we need introduce some new notations.
\begin{align*}
&V_{1}\left(\lambda_{0}\right)=
\begin{cases}
\frac{Y}{\lambda_{0}}-C^{*}\left(\lambda_{0}\right), & \lambda_{0} \in\left(0,\lambda_{max}\right)\\
\mathcal V\left(X\right), & \lambda_{0}=\lambda_{max}\quad \rm {where} \left\langle \mathcal V\left(X\right), C\right\rangle=\|X^{T}C\|_{2}.
\end{cases}\\
&V_{2}\left(\lambda,\lambda_{0}\right)=\frac{Y}{\lambda}-C^{*}\left(\lambda_{0}\right).\\
&V_{3}\left(\lambda,\lambda_{0}\right)=V_{2}\left(\lambda,\lambda_{0}\right)
-\frac{\left\langle V_{1}\left(\lambda_{0}\right),V_{2}\left(\lambda,\lambda_{0}\right)\right\rangle}{\|V_{1}\left(\lambda_{0}\right)\|^{2}_{F}}
V_{1}\left(\lambda_{0}\right).
\end{align*}
The idea of this process is same with Wang et al. (2015) in LASSO case, where the dual solution is a projection on a nonempty, closed and convex polytope. The dual solution of NMLR (1) is a projection on a set  that is not polytope. Therefore, the operator  $V_{1}\left(\lambda_{0}\right)$ doesn't have a closed form when $\lambda_{0}=\lambda_{max}$.
\subsection{PSRi}
This section will give a better tuning parameter selection rule than PSR based on the idempotence of the projection operator. We call it PSRi. Before showing the PSRi result, we need the following lemma, which gives a more accurate set containing $C^{*}(\lambda)$.
\begin{Lemma}
For a given $\lambda_{0}$, suppose the solution $C^{*}\left(\lambda_{0}\right)$ of  (3) is known. For any $ 0<\lambda< \lambda_{0}$,
 the dual solution $C^{*}\left(\lambda\right)$ can be estimated as follows
\begin{center}
$C^{*}\left(\lambda\right) \in \Omega_{1}\subseteq \Omega$,
\end{center}
where~~$\Omega_{1}:=\left\{C\Big|~\|C-C^{*}\left(\lambda_{0}\right)\|_{F}\leq \|V_{3}\left(\lambda,\lambda_{0}\right)\|_{F}\right\}$.
\end{Lemma}
In Lemma 3.1,  the radius of $\Omega_{1}$ is  $\rho=\|V_{3}\left(\lambda,\lambda_{0}\right)\|_{F}$ and
$\|V_{3}\left(\lambda,\lambda_{0}\right)\|_{F}\leq \left(\frac{1}{\lambda}-\frac{1}{\lambda_{0}}\right)\|Y\|_{F}$.
With the similar way of proving  Theorem 3.1, we get the PSRi theorem below.
\begin{Theorem}{\rm(\textbf{PSRi})}
Assume the solution $C^{*}\left(\lambda_{0}\right)$ of (3) is known for given $\lambda_{0}$.
For any $i\in \left\{1,\cdots,r\right\}$, if $0<\lambda <\lambda_{0}$ and
\begin{center}
$\sigma_{i}\left(X^{T}C^{*}\left(\lambda_{0}\right)\right)< 1-\|X\|_{2}\|V_{3}\left(\lambda,\lambda_{0}\right)\|_{F}$,
\end{center}
then the solution of NMLR (1) satisfies  rank$\left(B^{*}(\lambda)\right)\leq i-1$.
\end{Theorem}
As the similar arguments after Theorem 3.1, Theorem 3.3 needs the solution $C^{*}(\lambda_{0})$  is known and it may be difficult to solve for any given $\lambda_{0}$. Note that $C^{*}(\lambda_{max})=\frac{Y}{\lambda_{max}}$ from Proposition 3.1, so we have the following result.
\begin{Theorem}
 Suppose $X^{T}Y$ has at least two different singular values. For $i\in \{1,\cdots,r\}$, define $\lambda_{i}$ such that
\begin{center}
$\sigma_{i}\left(X^{T}Y\right)< \lambda_{max}-\lambda_{max}\|X\|_{2}\|V_{3}\left(\lambda_{i},\lambda_{max}\right)\|_{F}$.
\end{center}
Then the solution of NMLR  (1) satisfies rank$\left(B^{*}\left(\lambda_{i}\right)\right)\leq i-1$.
\end{Theorem}
From above results, it's sure that the performance of PSRi  in Theorems 3.3 and 3.4 is better than PSR. Because the dual solution of NMLR is a projection on a complex set, which is different from polytope in vector (see, Wang et al., 2015), $V_{3}\left(\lambda,\lambda_{max}\right)$ doesn't have a closed form. Hence, it's not easy to obtain a closed form of tuning parameters $\lambda_{1}, \cdots, \lambda_{r}$ as in Theorems 3.2.  Next, we give another tuning parameter selection rule based on a different property of the projection operator.
\subsection{PSRfn}
We get another tuning parameter selection rule PSRfn based on the firm nonexpansiveness of the projection operator in this section. In order to obtain the results, we first  give a lemma.
\begin{Lemma}
For a given $\lambda_{0}$, suppose the solution $C^{*}\left(\lambda_{0}\right)$ of  (3) is known.
For any $ 0<\lambda< \lambda_{0}$, the dual solution $C^{*}\left(\lambda\right)$ can be estimated as follows
\begin{center}
$C^{*}\left(\lambda\right) \in \Omega_{2}\subseteq \Omega$,
\end{center}
where~~$\Omega_{2}:=\left\{C\Big|~\left\|C-C^{*}\left(\lambda_{0}\right)
-\frac{1}{2}\left(\frac{1}{\lambda}-\frac{1}{\lambda_{0}}\right)Y\right\|_{F}
\leq \frac{1}{2}\left(\frac{1}{\lambda}-\frac{1}{\lambda_{0}}\right)\|Y\|_{F}\right\}$.
\end{Lemma}
In Lemma 3.2, the center of $\Omega_{2}$ is  $C^{*}\left(\lambda_{0}\right)+\frac{1}{2}\left(\frac{1}{\lambda}-\frac{1}{\lambda_{0}}\right)Y$ and the radius is  $\rho=\frac{1}{2}\left(\frac{1}{\lambda}-\frac{1}{\lambda_{0}}\right)\|Y\|_{F}$. By using the similar idea in Theorem 3.1, we get the PSRfn result.
\begin{Theorem}{\rm(\textbf{PSRfn})}
Assume  the solution  $C^{*}\left(\lambda_{0}\right)$ of (3) is known for given $\lambda_{0}$.
For any $i \in \{1,2,\cdots,r\}$, if $\lambda<\lambda_{0}$ and
\begin{center}
$\lambda>\frac{\lambda_{0}\|X\|_{2}\|Y\|_{F}}{2\lambda_{0}-2\lambda_{0}\sigma_{i}\left(X^{T}\left(C^{*}\left(\lambda_{0}\right)+\frac{1}{2}\left(\frac{1}{\lambda}-\frac{1}{\lambda_{0}}\right)Y\right)\right)+\|X\|_{2}\|Y\|_{F}}$,
\end{center}
then the solution of NMLR  (1) satisfies  $\sigma_{j}\left(B^{*}(\lambda)\right)=0$  $(j\geq i)$. Moreover, rank$\left(B^{*}(\lambda)\right)\leq i-1$.
\end{Theorem}
Theorem 3.5 needs that $C^{*}(\lambda_{0})$ is known, but it's not easy to compute for any $\lambda_{0}$. Next, we have a special result for $\lambda_{0}=\lambda_{max}$.
\begin{Theorem}
Suppose $X^{T}Y$ has  at least two different singular values. For $i\in \{1,\cdots,r\}$, let's define $\lambda_{i}$  as
\begin{center}
$\lambda_{i}=\frac{\|X^{T}Y\|_{2}\left(\|X\|_{2}\|Y\|_{F}+\sigma_{i}\left(X^{T}Y\right)\right)}
{2\|X^{T}Y\|_{2}-\sigma_{i}\left(X^{T}Y\right)+\|X\|_{2}\|Y\|_{F}}$.
\end{center}
Then, for any  $\lambda \in (\lambda_{i},\lambda_{i-1}]$ ($i\geq 2$), the solution of NMLR (1) satisfies
\begin{center}
 rank$(B^{*}(\lambda))\leq i-1$.
\end{center}
\end{Theorem}
The above results show that the performance of PSRfn in Theorems 3.5 and 3.6 is better than PSR.
\subsection{PSR+}
The results of PSRi and PSRfn are deduced separately from the idempotence and the firm nonexpansiveness properties of the projection operator. If one can combine the two properties together, a more accurate set that contains $C^{*}(\lambda)$ may be reached. So it does. We get the enhanced version PSR+ based on these two properties. Firstly, we give a lemma about the estimate set of $C^{*}\left(\lambda\right)$.
\begin{Lemma}
For a given $\lambda_{0}$, suppose the solution $C^{*}\left(\lambda_{0}\right)$ of  (3) is known.
For any $ 0<\lambda< \lambda_{0}$, the dual solution $C^{*}\left(\lambda\right)$ can be estimated as follows
\begin{center}
$C^{*}(\lambda)\in \Omega_{3}\subseteq \Omega$,
\end{center}
where~~$\Omega_{3}:=\left\{C\Big|\left\|C-C^{*}\left(\lambda_{0}\right)-\frac{1}{2}V_{3}\left(\lambda,\lambda_{0}\right)\right\|_{F}
\leq \frac{1}{2} \left\|V_{3}\left(\lambda,\lambda_{0}\right)\right\|_{F}\right\}$.
\end{Lemma}
In Lemma 3.3, we combine the  idempotence and the firm nonexpansiveness of projection operator to get an estimate set  $\Omega_{3}$. Its radius is $\rho=\frac{1}{2} \|V_{3}\left(\lambda,\lambda_{0}\right)\|_{F}$  and the center is $C^{*}\left(\lambda_{0}\right)
+\frac{1}{2}V_{3}\left(\lambda,\lambda_{0}\right)$.  For clarifying the relationship between $\lambda$ and $\rm rank(B^{*}(\lambda))$, we give the following theorem.
\begin{Theorem}{\rm(\textbf{PSR+})}
Assume  the solution of (3) $C^{*}\left(\lambda_{0}\right)$ is known for given $\lambda_{0}$.  For any $i\in \{1,\cdots,r\}$, if $0<\lambda<\lambda_{0}$ and
\begin{center}
$\sigma_{i}\left(X^{T}\left(C^{*}\left(\lambda_{0}\right)+\frac{1}{2}V_{3}\left(\lambda,\lambda_{0}\right)\right)\right)< 1-\frac{1}{2}\|X\|_{2}\|V_{3}\left(\lambda,\lambda_{0}\right)\|_{F}$,
\end{center}
then the solution of NMLR (1) satisfies  rank$\left(B^{*}(\lambda)\right)\leq i-1$ .
\end{Theorem}
Same as the analysis in Theorem 3.1, if the $C^{*}(\lambda_{0})$ can't be computed easily for any $\lambda_{0}$, we can choose $\lambda_{0}$ as $\lambda_{max}$. The result is presented in the next theorem.  We can't get a closed form of  $\{\lambda_{i}\}^{r}_{i=1}$ due to the speciality of $V_{3}\left(\lambda_{i},\lambda_{max}\right)$, but it doesn't influence the efficiency of this result.
\begin{Theorem}
Suppose $X^{T}Y$ has at least two different singular values. For $i\in \{1,\cdots,r\}$, define $\lambda_{i}$ such that
\begin{center}
$\sigma_{i}\left(X^{T}\left(\frac{Y}{\lambda_{max}}+\frac{1}{2}V_{3}\left(\lambda_{i},\lambda_{max}\right)\right)\right)< 1-\frac{1}{2}\|X\|_{2}\|V_{3}\left(\lambda_{i},\lambda_{max}\right)\|_{F}$.
\end{center}
Then the solution of NMLR (1) satisfies rank$\left(B^{*}\left(\lambda_{i}\right)\right)\leq i-1$.
\end{Theorem}
From Lemma 3.3,  we know that PSR+ outperform PSR. When $\lambda_{0}=\lambda_{max}$, the sequence of tuning parameters in PSR presents the closed form, while PSR+ doesn't.  The reason is the speciality of $V_{3}\left(\lambda_{i},\lambda_{max}\right)$, which is caused by the complex projection set.
\begin{Remark}
We claim the relationships among $\Omega, \Omega_{1}, \Omega_{2}$ and $\Omega_{3}$ as $\Omega_{3}\subseteq\Omega_{1}\subseteq\Omega$ and $\Omega_{3}\subseteq\Omega_{2}\subseteq\Omega$. From previous arguments, we have
$\Omega_{1}\subseteq \Omega$ and $\Omega_{2}\subseteq \Omega$. It remains to prove that $\Omega_{3}\subseteq \Omega_{1}$ and $\Omega_{3}\subseteq \Omega_{2}$. In fact, for any $C\in \Omega_{3}$, it means that
$\left\|C-C^{*}\left(\lambda_{0}\right)-\frac{1}{2}V_{3}\left(\lambda,\lambda_{0}\right)\right\|_{F}
\leq \frac{1}{2} \|V_{3}\left(\lambda,\lambda_{0}\right)\|_{F}.$ By Cauchy-Schwarz inequality,
$$\left\|C-C^{*}\left(\lambda_{0}\right)\right\|_{F}-\left\|\frac{1}{2}V_{3}\left(\lambda,\lambda_{0}\right)\right\|_{F}
\leq \left\|C-C^{*}\left(\lambda_{0}\right)-\frac{1}{2}V_{3}\left(\lambda,\lambda_{0}\right)\right\|_{F}
\leq \frac{1}{2} \|V_{3}\left(\lambda,\lambda_{0}\right)\|_{F}.$$
It implies that  $C\in \Omega_{1}$ which leads to $\Omega_{3}\subseteq \Omega_{1}$. Similarly, $\Omega_{3}\subseteq \Omega_{2}$. In  one and two dimensional  setting, these sets are showed in Figure 1. Hence, we get that PSRi, PSRfn and  PSR+ outperform  PSR, and PSR+ is the most accurate result. However, it's not clear which one of PSRi and PSRfn is better.
\begin{figure}[htbp]
\centering\includegraphics[height=4.2cm]{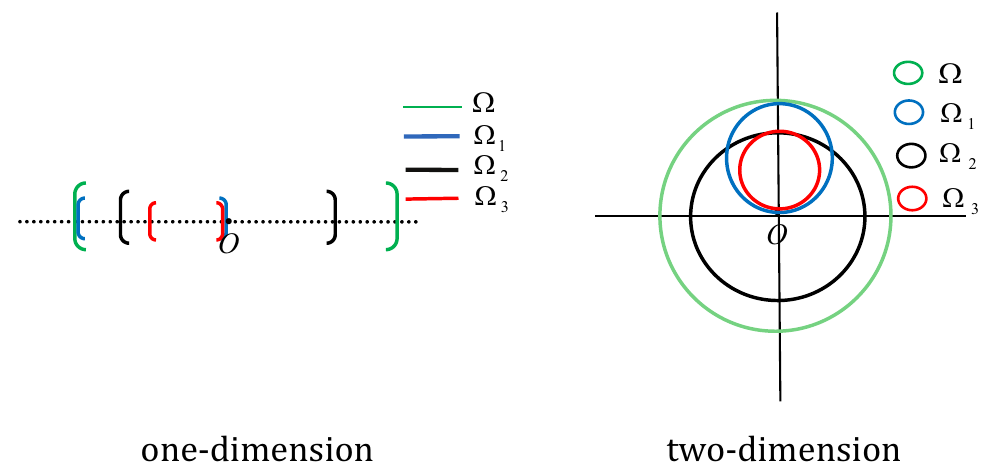}
\caption{The simple graphical representations}
\end{figure}
\end{Remark}
\section{Numerical studies}
In the previous section, we get the four tuning parameter selection rules for NMLR and these rules  depend on a solution of (3). Thus, an efficient algorithm for solving (3) is needed. Here, we present the popular first-order method, the alternating direction multiplier method (ADMM). See, e.g.,  Boyd et al. (2012), Fazek et al. (2013) and Bottou et al. (2018).  The numerical results on simulation and real data show that the four rules are all valuable and PSR+ is the most efficient one.

First, we give the detail process of ADMM for solving problem (3). We first transform (3) as a constraint problem
\begin{align*}
\underset{C,E} \min\left\{ \frac{1}{2}\|C-Y\|^{2}_{F}-\frac{1}{2}\|Y\|^{2}_{F}
+\delta_{\|\cdot\|_{2}\leq\lambda}\left(E\right)\right\}\quad s.t. \quad X^{T}C-E=0.
\end{align*}
Therefore, the augmented Lagrangian function is
$$\textit{L}_{\sigma}\left(C,E;Z\right)=\frac{1}{2}\|C-Y\|^{2}_{F}-\frac{1}{2}\|Y\|^{2}_{F}
+\delta_{\|\cdot\|_{2}\leq\lambda}\left(E\right)+\left\langle Z, X^{T}C-E\right\rangle
+\frac{\sigma}{2}\|X^{T}C-E\|^{2}_{F}.$$
We present the ADMM for (3) as follows.
\makeatletter\def\@captype{table}\makeatother
\begin{center}
\renewcommand\arraystretch{1}
\begin{tabular}{l}
\hline
\textbf{Algorithm:~~ADMM for solving problem (3)}\\
\hspace{2mm}
\textbf{Step 0:}\quad Set $C^{0},E^{0}$ and $Z^{0}$, let $\tau\in \left(0,\frac{1+\sqrt{5}}{2}\right)$ and $\sigma>0$;\\
\hspace{2mm}
\textbf{Step 1:}\quad Compute $C^{k+1}=\underset {C} {\rm{argmin}} \left\{\textit{L}_{\sigma}\left(C,E^{k};Z^{k}\right)\right\};$\\
\hspace{2mm}
\textbf{Step 2:}\quad Compute $E^{k+1}=\underset {E} {\rm{argmin}} \left\{\textit{L}_{\sigma}\left(C^{k+1},E;Z^{k}\right)\right\};$\\
\hspace{2mm}
\textbf{Step 3:}\quad Compute $Z^{k+1}=Z^{k}+\tau\sigma\left(X^{T}C^{k+1}-E^{k+1}\right).$\\
\hspace{2mm}
\textbf{Step 4:}\quad If a termination criterion is not met, go to Step 1-3.\\\hline
\end{tabular}
\end{center}
It's easy to get the closed form solutions for subproblems.
\begin{align*}
C^{k+1}
&=\underset {C} {\rm{argmin}}\left\{ \textit{L}_{\sigma}\left(C,E^{k};Z^{k}\right)\right\}
=\underset {C} {\rm{argmin}}  \left\{\frac{1}{2}\|C-Y\|^{2}_{F}+\left\langle X^{T}Z^{k},C\right\rangle\right\}
+\frac{\sigma}{2}\|X^{T}C-E^{k}\|^{2}_{F}\\
&=\left(I+\sigma XX^{T}\right)^{-1}\left(\sigma XE^{k}+Y-XZ^{k}\right).\\
E^{k+1}
&=\underset {E} {\rm{argmin}}\left\{\textit{L}_{\sigma}\left(C^{k+1},E;Z^{k}\right)\right\}
=\underset {E} {\rm{argmin}}\left\{\delta_{\|\cdot\|_{2}\leq\lambda}\left(E\right)
+\frac{\sigma}{2}\left\|E-X^{T}C^{k+1}-\frac{Z^{k}}{\sigma}\right\|^{2}_{F}\right\}\\
&=\Pi_{\|\cdot\|_{2}\leq \lambda}\left(X^{T}C^{k+1}+\frac{Z^{k}}{\sigma}\right).
\end{align*}
The convergence of two-blocks ADMM  is well-known. For the special case (3), we describe its convergence result as follows.
\begin{Theorem}
Assume that the solution set of (3) is nonempty. Let $\left\{(C^{k},E^{k},Z^{k})\right\}$ be generated from ADMM for $\tau \in (0,\frac{1+\sqrt{5}}{2})$. Then the sequence $\left\{(C^{k},E^{k})\right\}$ converges to the solution of problem (3) and $\left\{Z^{k}\right\}$ converges to the solution of problem (2).
\end{Theorem}
\subsection{Simulation}
We evaluate tuning parameter selection rules on simulation data. We randomly simulate matrix X and B distributed as standard norm distribution.  The dimension of X and B are designed as $X\in R^{100\times 5000}$ and $B\in R^{5000\times 500}$. It means that the sample size is 100, the dimensions of prediction and response variables are 5000 and 500, respectively. Each column of random error $W$ has mean 0 and standard variance 0.01. According to $Y=XB+W$,  the response matrix is gotten. Because $V_{3}\left(\lambda,\lambda_{max}\right)$ has no closed form, we omit the results of $\lambda_{0}=\lambda_{max}$.
We present the results of $\lambda_{0}=0.1\lambda_{max}$. Denoting $\lambda$, $\lambda_{i}$, $\lambda_{fn}$ and $\lambda_{+}$ as the tuning parameters under PSR, PSRi, PSRfn and PSR+, respectively.

All results of them are presented in Figure 2 and Table 1.  Figure 2 gives a simple show of the four tuning parameter selection rules, the measurement is the rank of $B^{*}(\lambda)$.  Although Figure 2  shows a nonincreasing trend of the rank of $B^{*}(\lambda)$ with $\lambda$ increasing, a certain rank corresponds to different tuning parameter under these rules where $\lambda_{+}$  is the smallest and $\lambda$ largest. In order to clearly present this results, Figure 2 also gives the rank of $B^{*}(\lambda)$ ranging from 30 to 50. Table 1 gives numerical values of rank and tuning parameters, they are accord with Figure 2. Therefore, we prove our tuning parameter selection rules are valuable and PSR+ performs best among these rules.
\begin{figure}[htbp]
\caption{The rank of B with different $\lambda$}
\centering
\subfigure[All results]{
\centering
\includegraphics[width=3.3in]{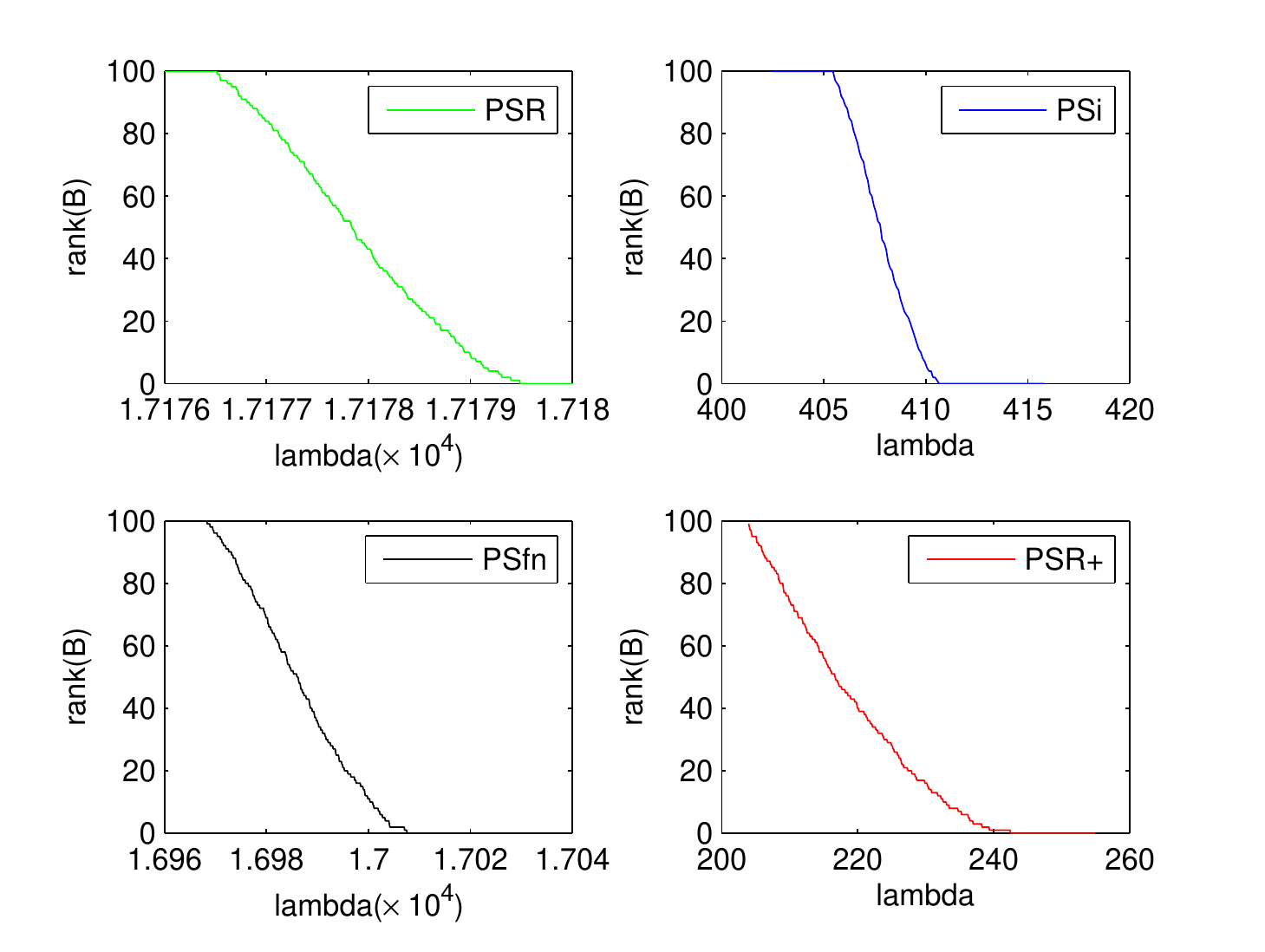}
}

\subfigure[Part results]{
\centering
\includegraphics[width=3.3in]{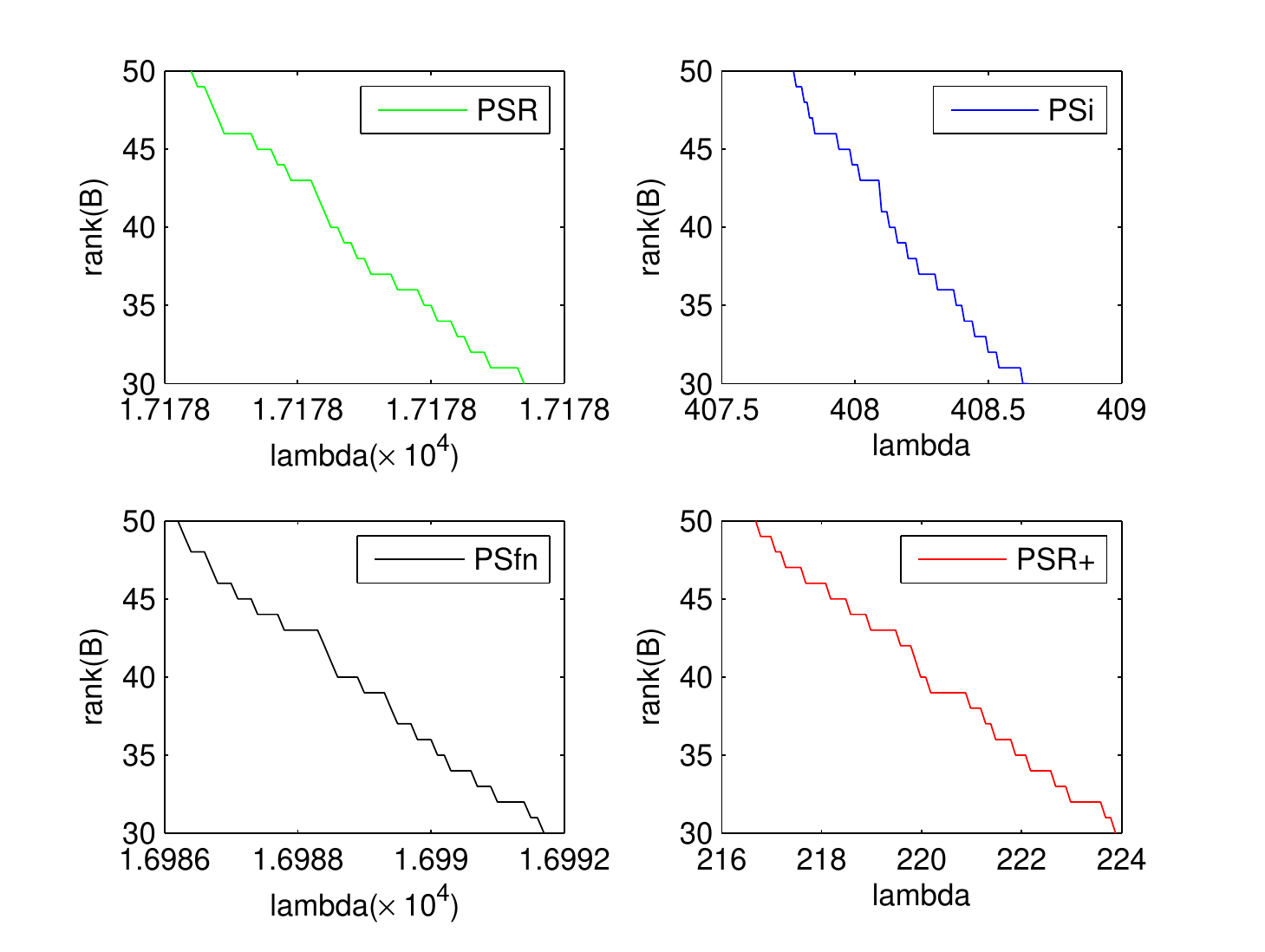}
}
\end{figure}
\makeatletter\def\@captype{table}\makeatother
\begin{center}
\renewcommand\arraystretch{0.8}
\begin{tabular}{|c|c|c|c|c|c|c|c|c|c|}
\hline
$\rm rank(B^{*}(\lambda))$&100&50&25&13&6&3&2&1&0\\\hline
$\lambda$ (17176+)&0.51&1.84&2.47&3.90&3.11&3.28&3.31&3.40&3.49\\\hline
$\lambda_{i}$ (405+)&0.46&2.16&3.86&3.96&4.46&5.26&5.36&5.56&5.57\\\hline
$\lambda_{fn}$ (16968+)&0.30&18.00&25.70&31.30&34.10&36.10&36.2&39.1&39.6\\\hline
$\lambda_{+}$ (203+)&0.98&13.68&22.88&27.88&32.29&33.98&35.28&36.38&39.48\\\hline
\end{tabular}
\caption{The numerical values of Figure 2}
\end{center}
\subsection{Real data}
In this section, we evaluate the algorithm for solving (3) and tuning parameter selection rules on a picture dataset that contains different shape black-and-white pictures.

First, the picture information are input as B. Then, we simulation X satisfying that each column distributes the standard norm distribution. The error matrix $W$ is simulated as norm distribution with mean 0 and standard variance 0.01. According to $Y=XB+W$, the response matrix is obtained. The picture recovery results are presented in Figure 3 where each subfigure includes the real picture in the left and recovery picture right. In Table 2, we report the dimensions of these pictures and some measurements for evaluating the algorithm, including time, iteration and MSE defined as MSE$=\frac{\|B-Z^{k}\|_{F}^{2}}{pq}$.

\makeatletter\def\@captype{figure}\makeatother
\begin{figure}[htbp]
\caption{The comparison results}
\centering
\subfigure[device0-14]{
\begin{minipage}[c]{0.2\textwidth}
\centering
\includegraphics[width=1in]{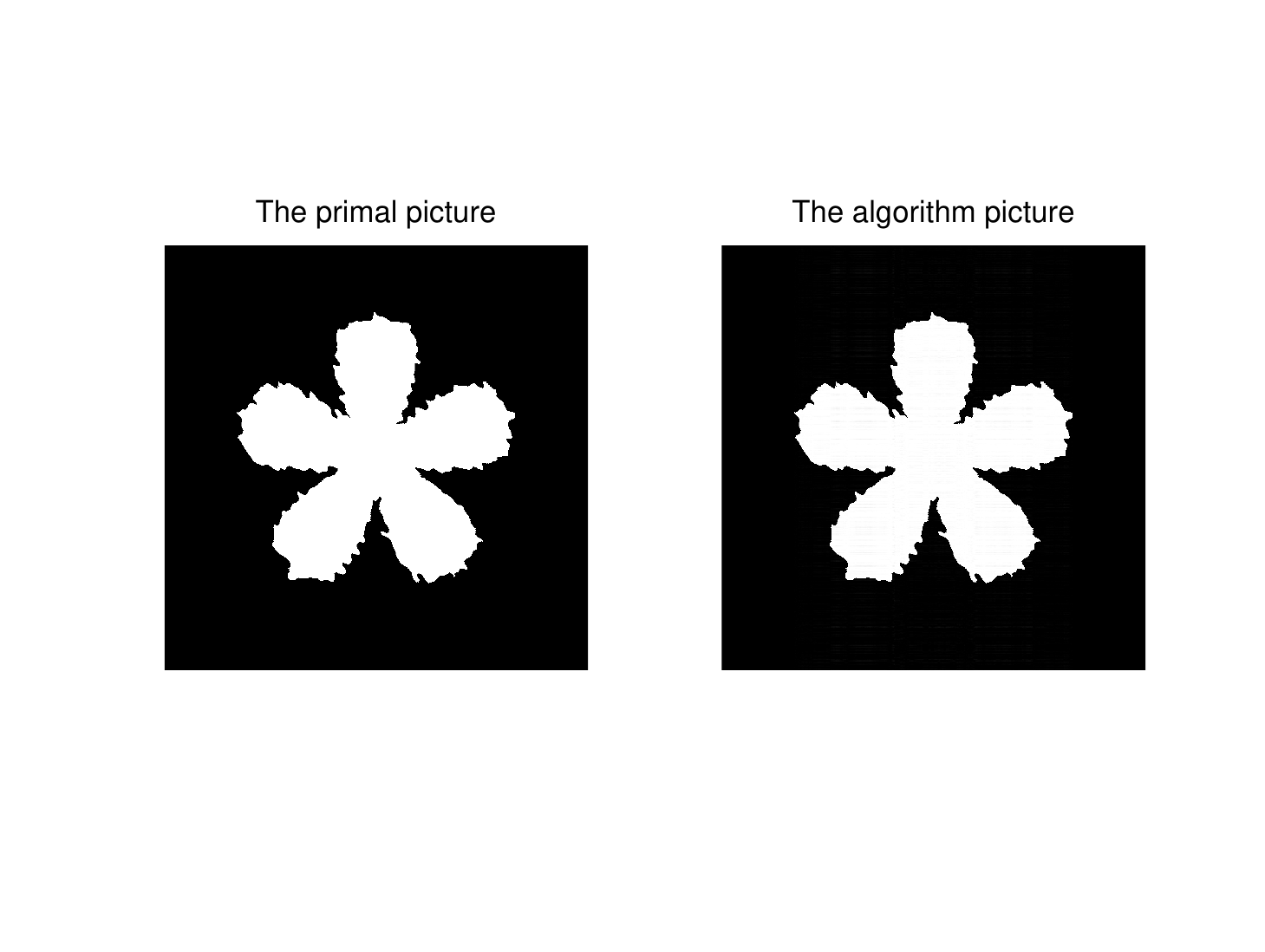}
\end{minipage}
}
\subfigure[fly-8]{
\begin{minipage}[c]{0.2\textwidth}
\centering
\includegraphics[width=1in]{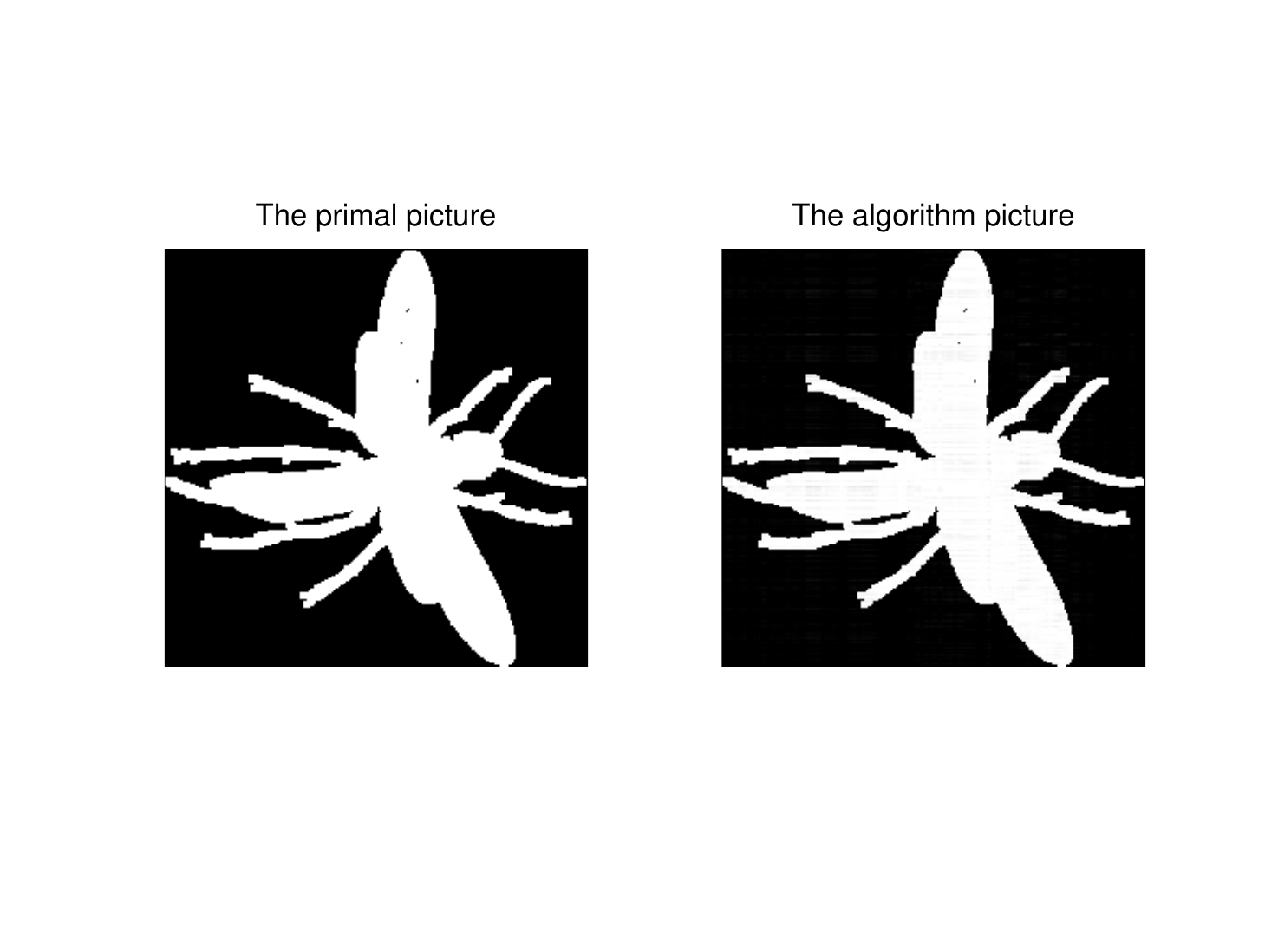}
\end{minipage}
}
\subfigure[butterfly-10]{
\begin{minipage}[c]{0.2\textwidth}
\centering
\includegraphics[width=1in]{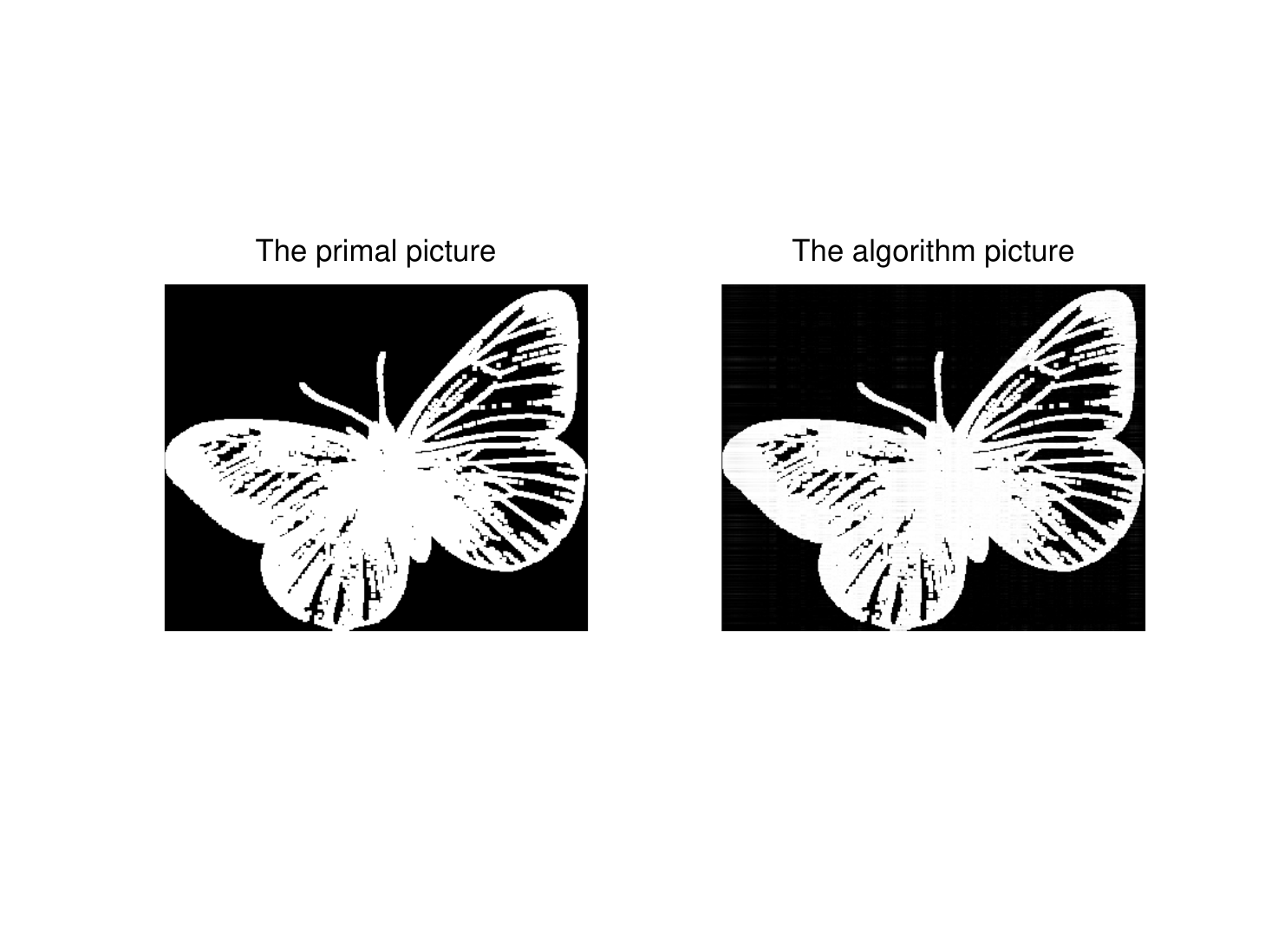}
\end{minipage}
}
\subfigure[turtle-14]{
\begin{minipage}[c]{0.2\textwidth}
\centering
\includegraphics[width=1in]{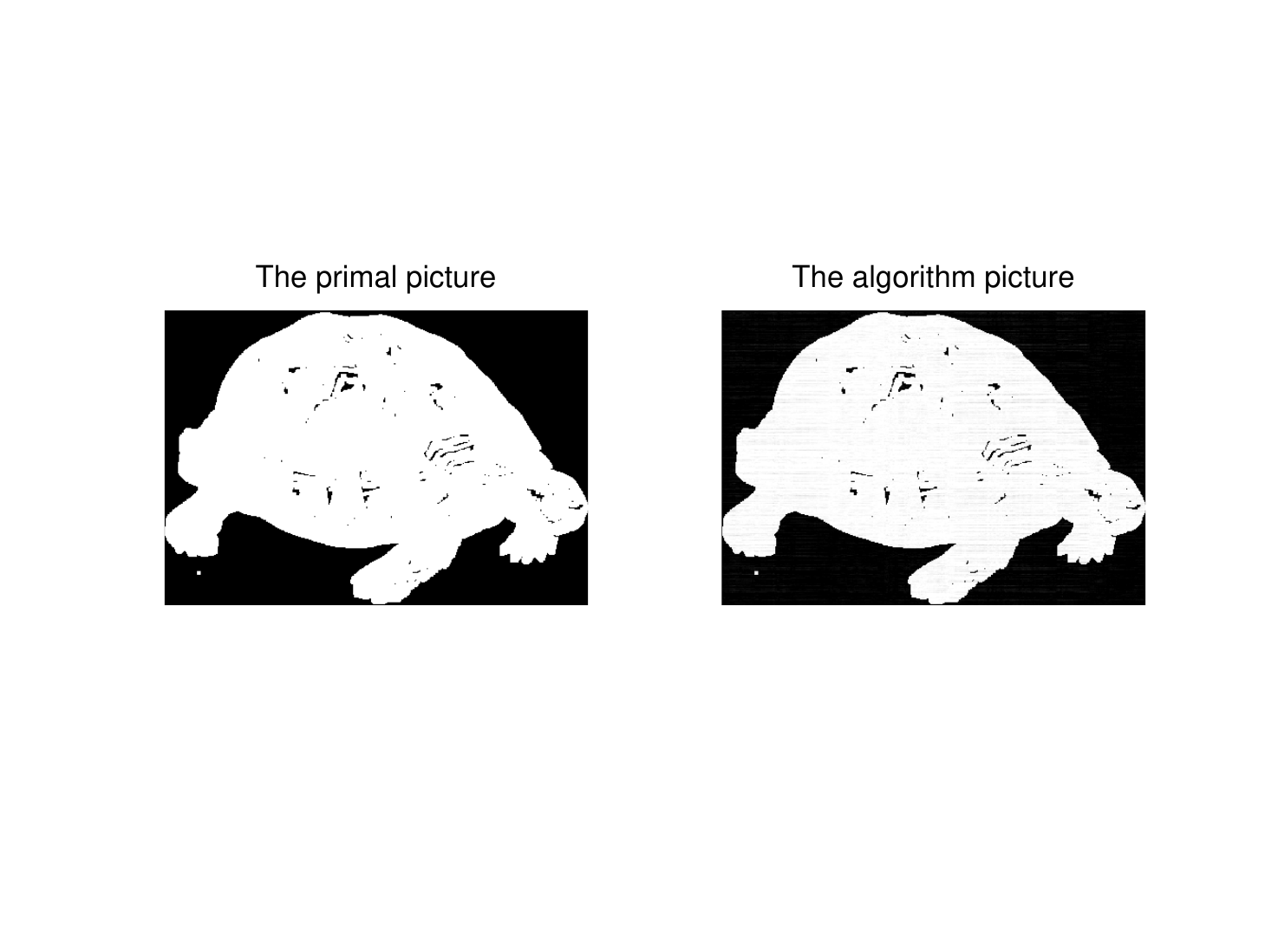}
\end{minipage}
}

\subfigure[bat-4]{
\begin{minipage}[c]{0.2\textwidth} 
\centering
\includegraphics[width=1in]{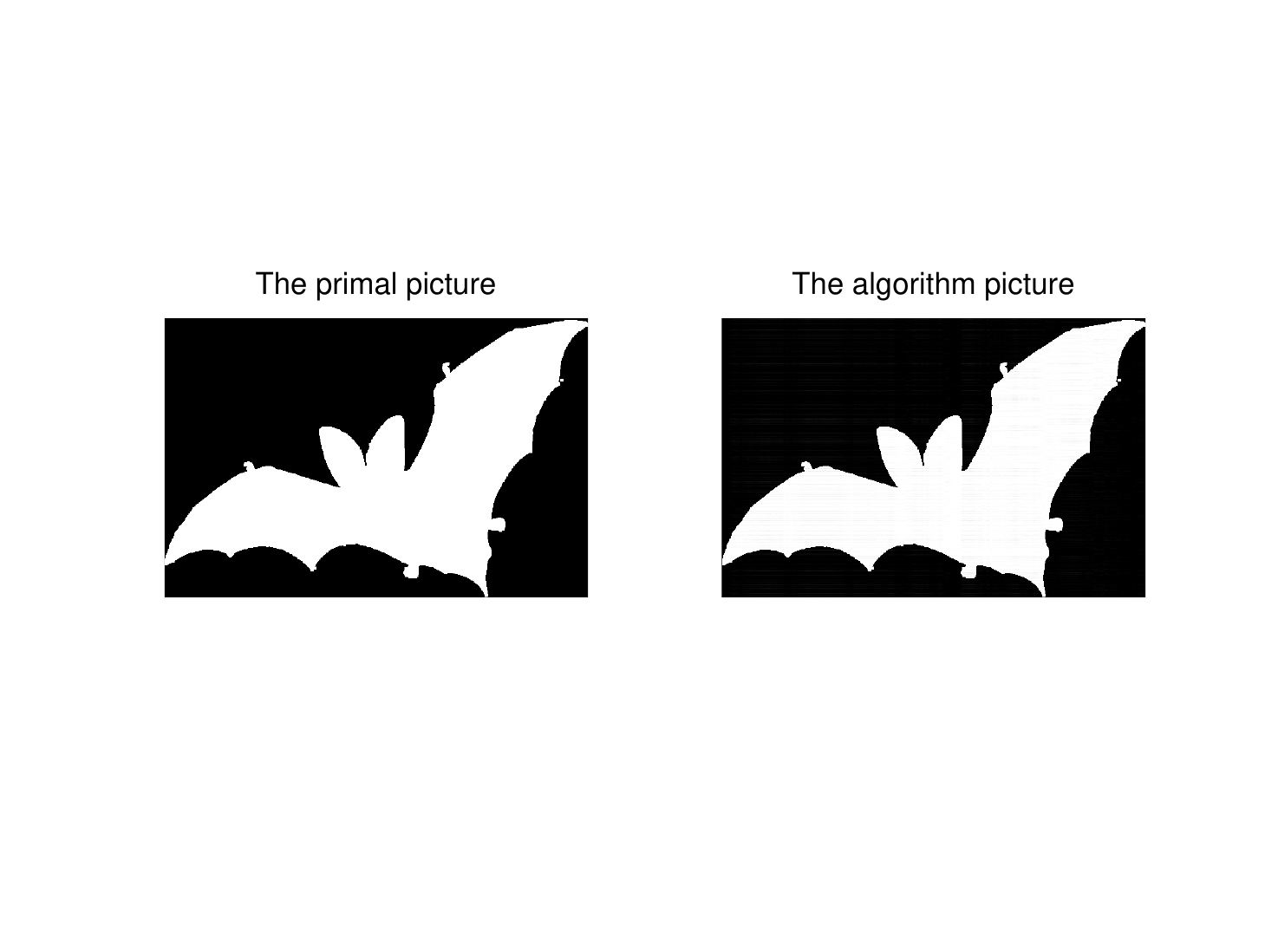} 
\end{minipage}
}
\subfigure[hat-10]{
\begin{minipage}[c]{0.2\textwidth}
\centering
\includegraphics[width=1in]{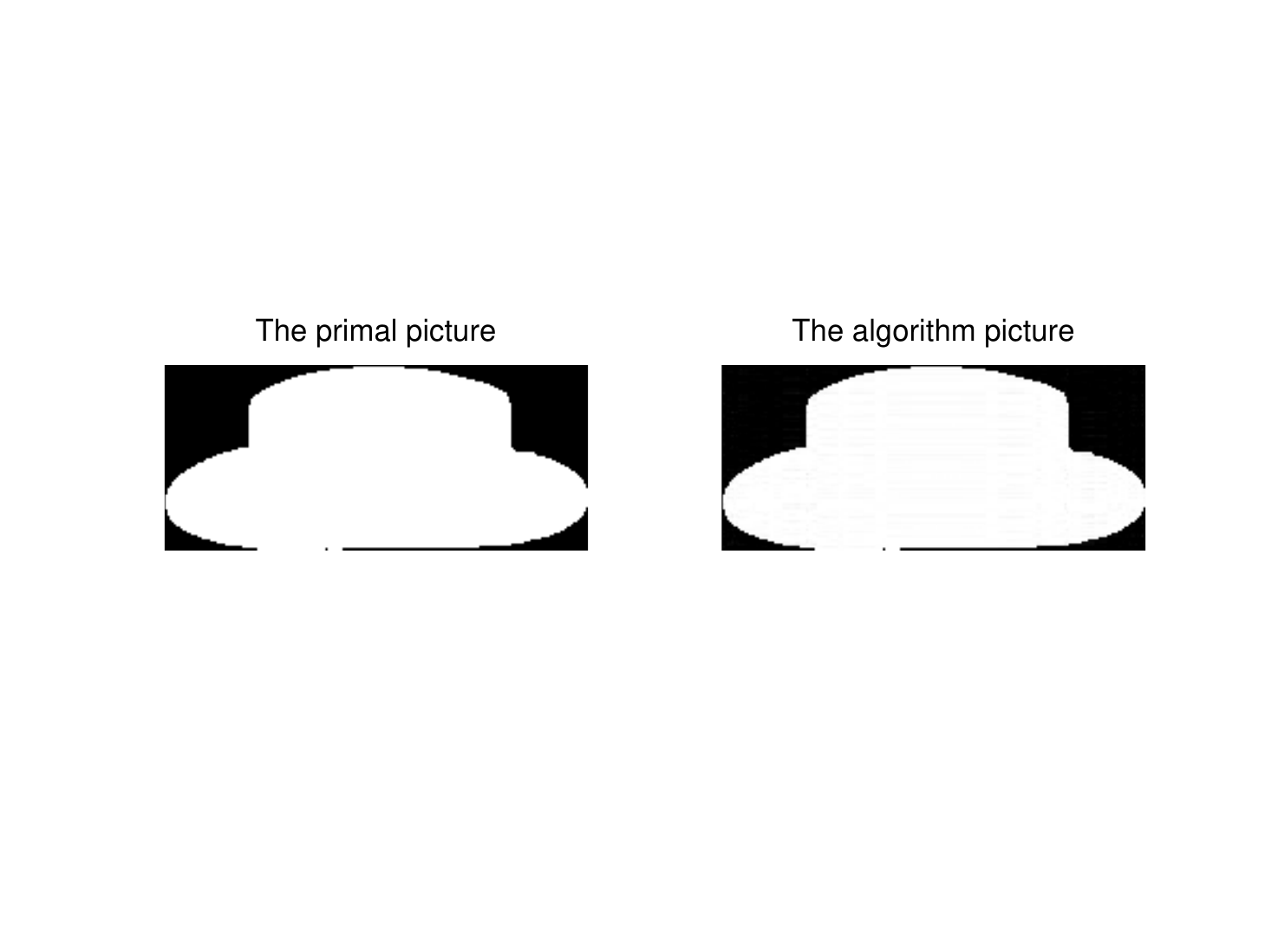}
\end{minipage}
}
\subfigure[lizzard-3]{
\begin{minipage}[c]{0.2\textwidth}
\centering
\includegraphics[width=1in]{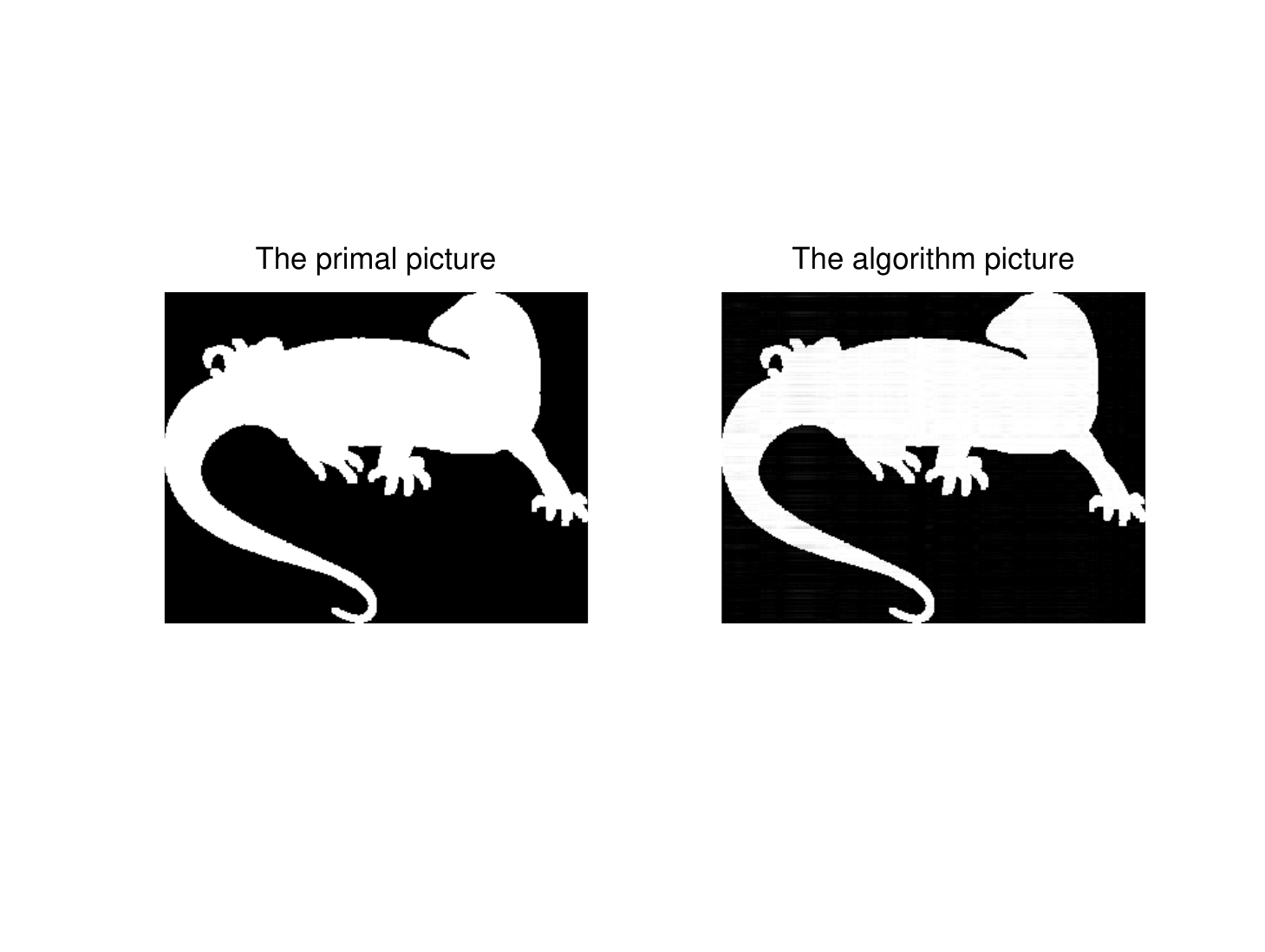}
\end{minipage}
}
\subfigure[pocket-20]{
\begin{minipage}[c]{0.2\textwidth}
\centering
\includegraphics[width=1in]{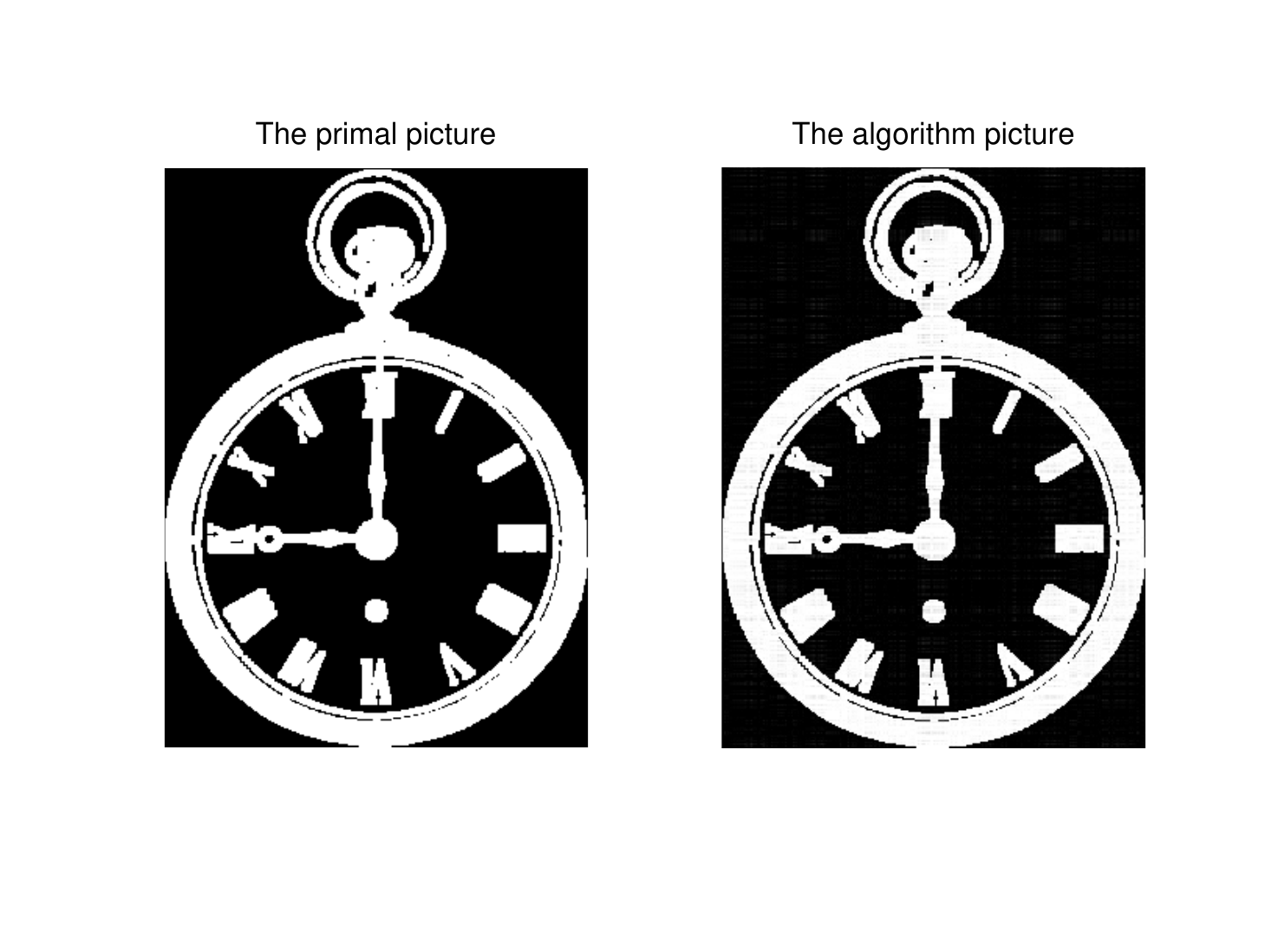}
\end{minipage}
}
\end{figure}
\makeatletter\def\@captype{table}\makeatother
\begin{center}
\renewcommand\arraystretch{0.8}
\begin{tabular}{|c|c|c|c|c|}
\hline
name&time(s)&Iterations&MSE &$\lambda$\\\hline
device0-14&88.9195&78&7.7141e-005&0.1604\\\hline
fly-8&7.1075&107&3.6480e-004&0.0295\\\hline
butterfly-10&20.4577&103&0.0016&0.6989\\\hline
turtle-14&38.5631&153&9.4783e-004&1.1074\\\hline
bat-4&61.7779&85& 1.8973e-004&0.1424\\\hline
hat-10&1.0372&94& 9.6308e-005&0.1354\\\hline
lizzard-3&11.4488&92&3.3952e-005&0.0563\\\hline
pocket-20&13.1728&77&2.8740e-004&0.0531\\\hline
\end{tabular}
\caption{The picture recovery Results}\label{table}
\end{center}
\noindent Figure 3 and Table 2 demostrate that ADMM  is efficient to solve dual problem (3). Therefore,  this algorithm can be used to provide the dual solution in tuning parameter selection rules. Then, the tuning parameter selection rules in Section 3 can be verified on this real picture dataset. Next, we show the performances of PSR, PSRi, PSRfn and PSR+ on picture device0-14 and pocket-20. In order to do so, we choose $\lambda_{0}=0.5\lambda_{max}$. For each picture, there are two results, which include the numerical values and figure of tuning parameters and the corresponding rank. The tuning parameter selection rules are proved valuable on these two pictures and PSR+ is the most efficient one.
\makeatletter\def\@captype{table}\makeatother
\begin{center}
\renewcommand\arraystretch{0.8}
\begin{tabular}{|c|c|c|c|c|c|c|c|c|c|c|c|}
\hline
$\rm r(B^{*}(\lambda))$&512&256&128&64&32&16&8&4&2&1&0\\\hline
$\lambda$ (55770+)&8&9&14&20&31&48&67&118&320&363&1021\\\hline
$\lambda_{i}$ (14660+)&1&2&7&13&23&39&56&101&290&349&967\\\hline
$\lambda_{fn}$ (45330+)&3&9&60&138&270&438&723&1194&3180&3645&8349\\\hline
$\lambda_{+}$ (8155+)&0.7&1.7& 11.7&23.7&43.7&65.7&112.7&203.7&586.7&665.7&1770.7\\\hline
\end{tabular}
\caption{The numerical values of device0-14 $(512\times 512)$}
\end{center}
\begin{figure}[htbp]
\caption{The rank of B under different $\lambda$ on device0-14}
\centering
\subfigure[All results]{
\centering
\includegraphics[width=4in]{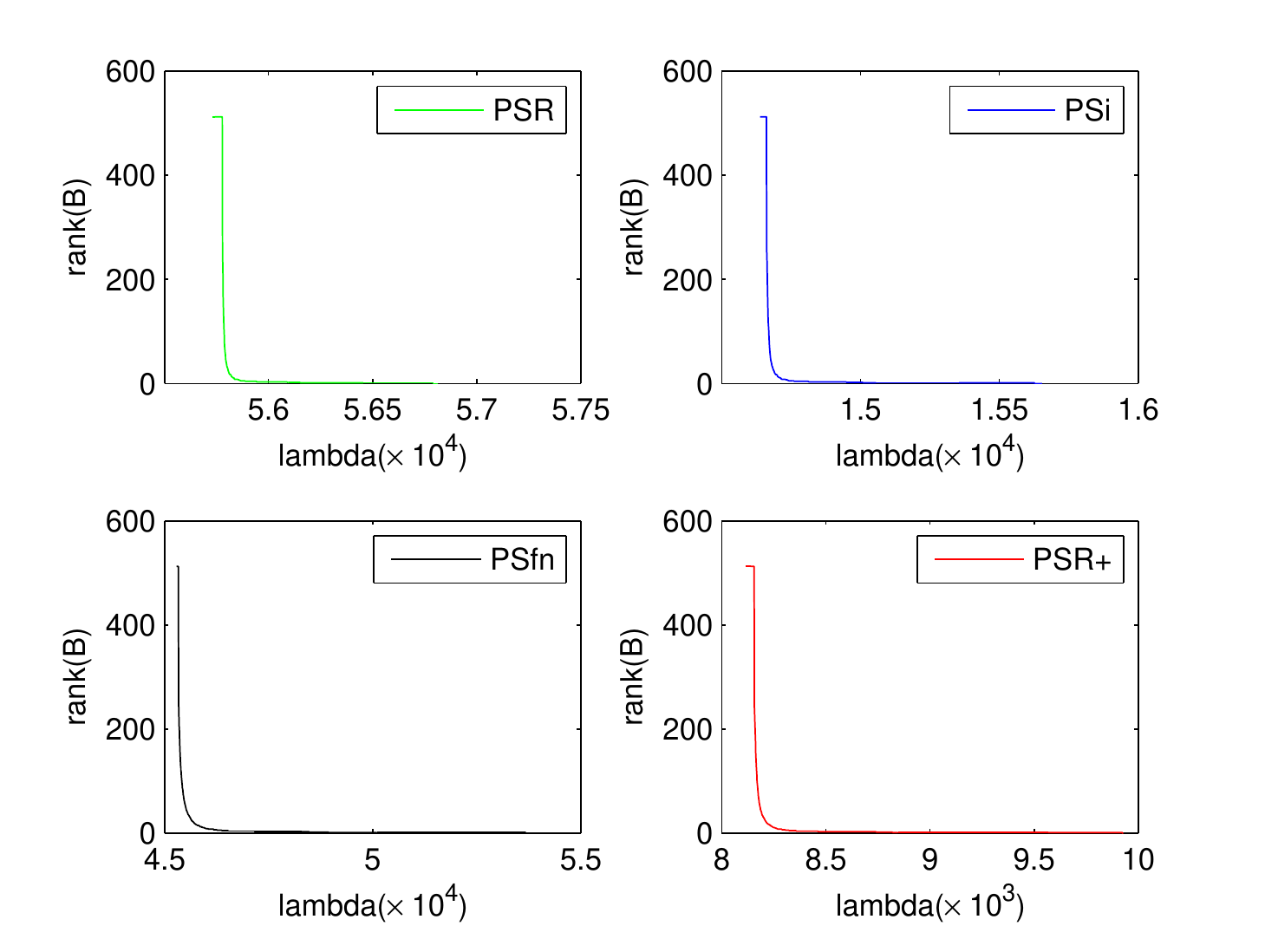}
}

\subfigure[Part results]{
\centering
\includegraphics[width=4in]{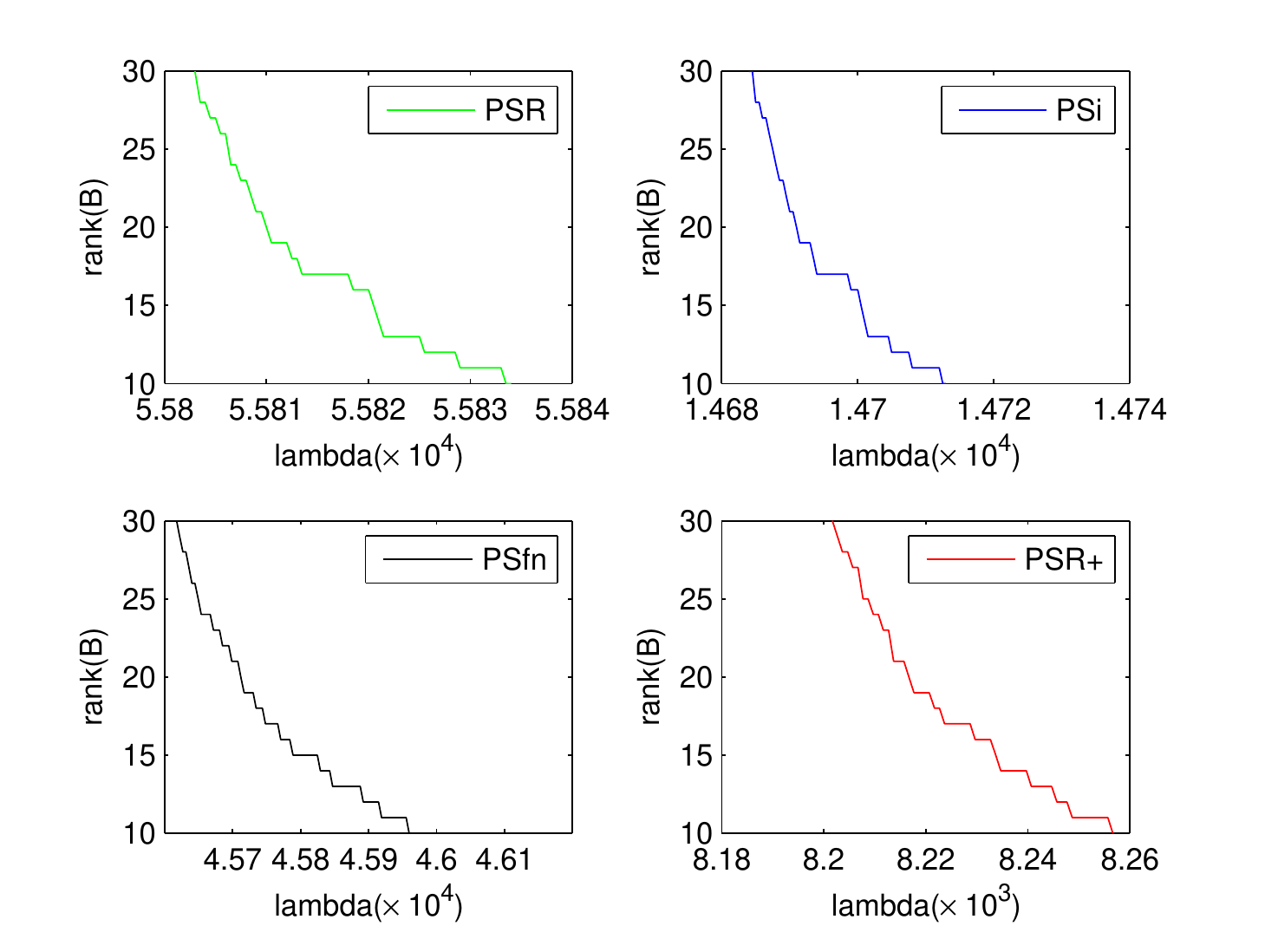}
}
\end{figure}
\makeatletter\def\@captype{table}\makeatother
\begin{center}
\renewcommand\arraystretch{0.8}
\begin{tabular}{|c|c|c|c|c|c|c|c|c|c|c|}
\hline
$\rm rank(B^{*}(\lambda))$&272&136&68&34&17&8&4&2&1&0\\\hline
$\lambda$ (24500+)&15&22&30&41&57&87&64&236&323&460\\\hline
$\lambda_{i} (8274+)$&0.8&9.8&19.8&33.8&55.8&94.8&197.8&292.8&414.8&612.8\\\hline
$\lambda_{fn}$ (20600+)&14&56&120&206&338&500&1090&1476&2052&2279\\\hline
$\lambda_{+}$ (4792+)&0.9&13.9&29.9&50.9&85.9&142.9&290.9&442.9&608.9&845.9\\\hline
\end{tabular}
\caption{The numerical values of pocket-20 $(372\times 272)$}
\end{center}
\begin{figure}[htbp]
\caption{The rank of B under different $\lambda$ on pocket-20}
\centering
\subfigure[All results]{
\centering
\includegraphics[width=4in]{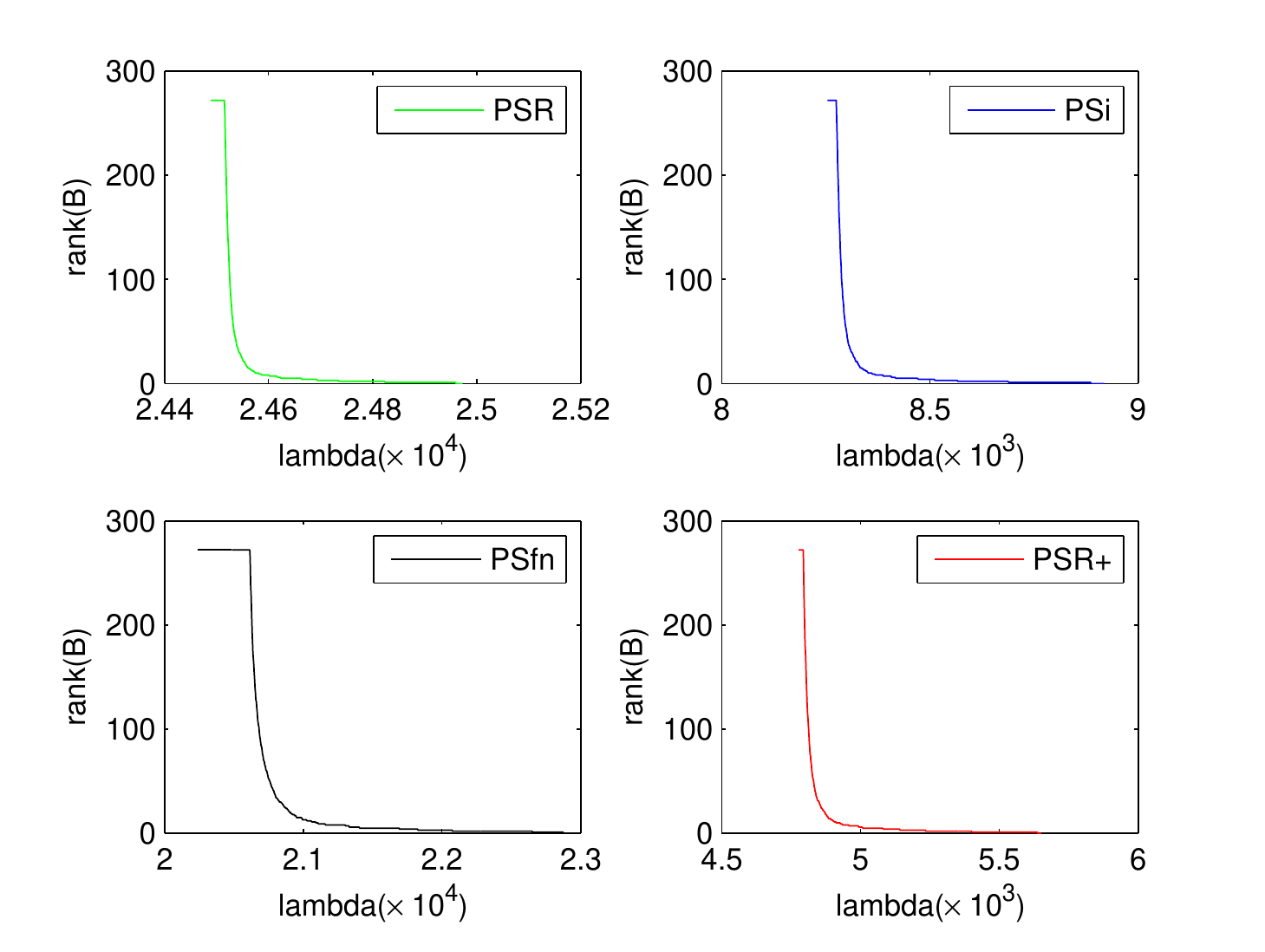}
}

\subfigure[Part results]{
\centering
\includegraphics[width=4in]{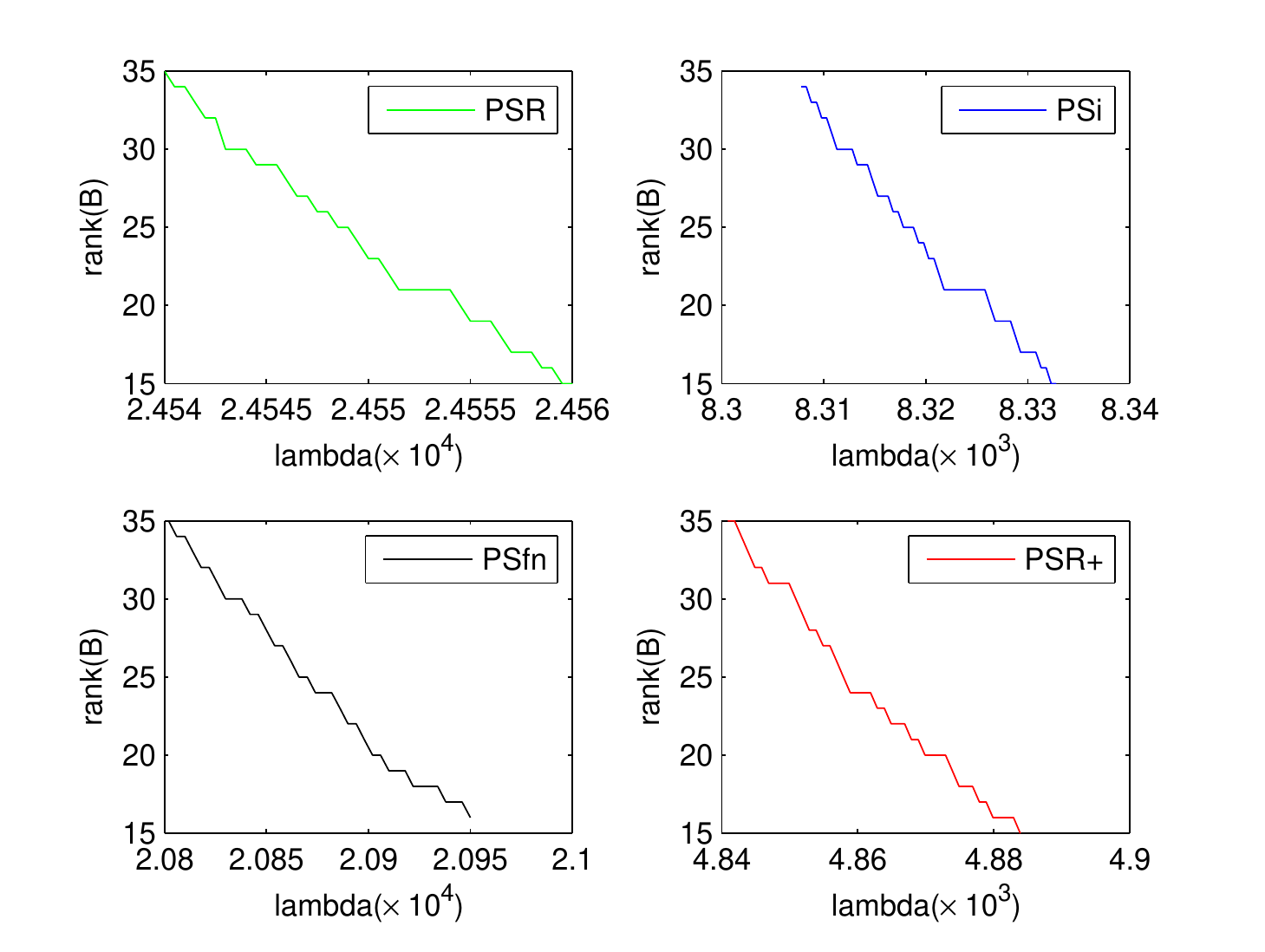}
}
\end{figure}
\section{Conclusion}
With the help of optimization techniques, this paper focus on the tuning parameter selection rules for nuclear norm regularized  multivariate linear regression (NMLR) in high-dimensional setting. We claim that the tuning parameter selection is closely related to the dual solution of NMLR. Then, we build four rules PSR, PSRi, PSRfn and PSR+, and discuss about the relationships among them. Moreover, we give a sequence of tuning parameters and the corresponding intervals, which states  that the rank of the estimation coefficient matrix is no more than a fixed number for the tuning parameter in the given interval.  Furthermore, we design an efficient ADMM to solve the dual problem of NMLR and our rules are illustrated to be valuable on simulation and real data. Actually, our rules are applicable to any efficient algorithm for solving the dual problem of NMLR.
\section*{Acknowledgements}
This work was supported by the National Science Foundation of China (11431002, 11671029). 
\section*{Appendix A: Main concepts and the duality theory}
\subsection*{Appendix A1: The concepts of conjugate function and projection operator}
The following definitions and results are  from Rockafellar (1970).
\begin{Definition}
Let $f: R^{p\times q}\rightarrow R$, the conjugate function $f^{*}: R^{p\times q}\rightarrow R$ of $f$ is defined as
$$f^{*}(M)= \underset {N\in dom(f)} \sup\left\{\langle M,N\rangle-f(N)\right\}.$$
\end{Definition}
\noindent If $f(M)=\|M\|_{*}$, we can get $f^{*}(M)= \underset {N \in R^{p\times q}} \sup\left\{\langle M,N\rangle-\|N\|_{*}\right\}=\delta_{\|\cdot\|_{2}\leq 1}(M)$, where $\delta_{\|\cdot\|_{2}\leq 1}(M)$ is an indicator function  defined as
$$\delta_{\|\cdot\|_{2}\leq 1}(M)=
\begin{cases} 0,&\|M\|_{2}\leq 1\\
               +\infty, & \|M\|_{2}> 1.\\
\end{cases}$$
\noindent If $f(M)=\|M\|_{F}$,  $f^{*}(M)= \underset {N\in R^{p\times q}} \sup\left\{\langle M,N\rangle-\|N\|_{F}\right\}=\delta_{\|\cdot\|_{2}\leq 1}(M)$.
\begin{Definition}
For an arbitrary vector $\omega$ and a convex set $\Omega$, the projection operator $P_{\Omega}(\cdot)$ is defined as
\begin{center}
$P_{\Omega}(\omega)=\underset{\mu \in \Omega}{\rm{argmin}}~~\|\mu-\omega\|^{2}_{F}$.
\end{center}
\end{Definition}
The following lemma gives an equivalent definition of the projection operator.
\begin{Lemma}
Suppose $\Omega$ is a nonempty, closed and convex set, $\mu=P_{\Omega}\left(\omega\right)$ if and only if
\begin{center}
$\left\langle \omega-\mu,\upsilon-\mu\right\rangle \leq 0$~~~~~~for any $\upsilon \in \Omega$.
\end{center}
\end{Lemma}
\noindent Here are some basic properties of projection operator.
\begin{Lemma}
Let $\Omega$ be any nonempty, closed and convex set, then the projection operator on $\Omega$ is \\
(1) nonexpansive, i.e.,
$\|P_{\Omega}(\omega_{2})-P_{\Omega}(\omega_{1})\|\leq \|\omega_{2}-\omega_{1}\|
~~~\forall \omega_{2},\omega_{1} \in \Omega.$\\
(2) idempotent, i.e.,
$P_{\Omega}(P_{\Omega}(\omega))=P_{\Omega}(\omega)~~~\forall \omega \in \Omega.$\\
(3) firmly nonexpansive, i.e.
$\|P_{\Omega}\left(\omega_{1}\right)-P_{\Omega}\left(\omega_{2}\right)\|^{2}+
\|\left(I-P_{\Omega}\right)\left(\omega_{1}\right)-\left(I-P_{\Omega}\right)\left(\omega_{2}\right)\|^{2}
\leq \|\omega_{1}-\omega_{2}\|^{2}\forall \omega_{2},\omega_{1} \in \Omega,$
where $I$ is the identity operator.
\end{Lemma}
\noindent There is another property of the projection operator which is showed in Lemma 5.3 below.
\begin{Lemma}
Let $\Omega$ be any convex set. $\nu=P_{\Omega}\left(\omega\right)+t\left(\omega-P_{\Omega}\left(\omega\right)\right)$ for $\forall \omega \in \Omega$ and $t\geq 0$. It holds that
\begin{center}
$P_{\Omega}\left(\nu\right)=P_{\Omega}\left(\omega\right)$.
\end{center}
\end{Lemma}
\subsection*{Appendix A2: The dual theory of NMLR}
First, we show the dual problem of (2). The Lagrangian function of (2) is
\begin{center}
$\textit{L}\left(B,A;\widetilde{C}\right)=\lambda||B||_{*}+\frac{1}{2}\|A\|^{2}_{F}+\left\langle \widetilde{C}, Y-XB-A\right\rangle$.
\end{center}
where $\widetilde{C}\in R^{n\times q}$ is a Lagrangian multiplier. We have the Lagrangian dual problem of (2)
$$\underset{\widetilde{C}}\max~\underset{B,A}\min\left\{\textit{L}\left(B,A;\widetilde{C}\right)\right\}.$$
It's not hard to yield the closed form of $\underset{B,A}\min\textit{L}\left(B,A;\widetilde{C}\right)$ as follows.
\begin{align*}
\underset{B,A}\min\textit{L}\left(B,A;\widetilde{C}\right)&=\underset{B,A}\min \left\{\lambda||B||_{*}+\frac{1}{2}\|A\|^{2}_{F}+\left\langle \widetilde{C}, Y-XB-A\right\rangle\right\}\\
&=\underset{B}\min\left\{\lambda||B||_{*}-\left\langle X^{T}\widetilde{C}, B\right\rangle\right\}+
\underset{A}\min \left\{\frac{1}{2}\|A\|^{2}_{F}-\left\langle \widetilde{C}, A\right\rangle\right\}+\left\langle \widetilde{C}, Y\right\rangle\\
&=-\delta_{\|\cdot\|_{2}\leq\lambda}\left(X^{T}\widetilde{C}\right)-\frac{1}{2}\|\widetilde{C}\|^{2}_{F}+\left\langle \widetilde{C}, Y\right\rangle.
\end{align*}
The last equality is a direct result of conjugate function. Thus the dual problem of (2) is
\begin{align*}
\underset{\widetilde{C}}\max \left\{-\frac{1}{2}\|\widetilde{C}-Y\|^{2}_{F}+\frac{1}{2}\|Y\|^{2}_{F}\right\}
\quad s.t. \quad \|X^{T}\widetilde{C}\|_{2}\leq\lambda.
\end{align*}
Taking $C=\frac{\widetilde{C}}{\lambda}$, we have
\begin{align*}
\underset{C}\min \left\{ \frac{\lambda^{2}}{2}\left\|C-\frac{Y}{\lambda}\right\|^{2}_{F}-\frac{1}{2}\|Y\|^{2}_{F}\right\}
\quad s.t. \quad \|X^{T}C\|_{2}\leq1,
\end{align*}
which is the dual problem of (2).

\textbf{Proof of Theorem 2.1}
\begin{proof}
Now we discuss the relationship between the convex optimization problem (2) and its dual (3). The objective function of (2) is $f:=\left\{\lambda||B||_{*}+\frac{1}{2}\|A\|^{2}_{F}\right\}$ and the feasible area $S:=\left\{(B,A)\Big|Y-XB-A=0\right\}$. For convex optimization problems with linear constraints, there is an important assumption named Slater's constraint qualification. If a convex optimization problem satisfies Slater's CQ, it follows from  Rockafellar (1970) that  the solutions of primal and dual problems are KKT points.

\textbf{Slater's CQ}: There exists $\theta \in$ ri(dom $(f)) \bigcap S$, where $f$ is the objective function and S is the feasible area of optimization problem.

\indent It's clear that there exist $B=0,A=Y$ such that $Y-XB-A=0$, which means (2) satisfies Slater's CQ. Because $\Omega_{D}$ is a nonempty, closed and convex set and $C^{*}(\lambda)=P_{\Omega_{D}}(\frac{Y}{\lambda})$, it's sure that problem (3)  have a solution. By solving (4) under $C=C^{*}(\lambda)$, we obtain the solution of (2). So, based on the Rockafellar (1970), the strong duality theorem holds on problems (2) and (3).
\end{proof}
\section*{Appendix B: The proofs of results in Section 3}
\textbf{Proof of Proposition 3.1}
\begin{proof}
We first prove the "only if"  part. Based on the KKT system (5),  it's obvious that if $B^{*}(\lambda)=0$, the solution of problem (3) is $$C^{*}(\lambda)=\frac{Y}{\lambda}.$$ It means $\frac{Y}{\lambda}\in \Omega_{D} $, which implies $\|X^{T}\frac{Y}{\lambda}\|_{2}\leq 1$. That is
$$\|X^{T}Y\|_{2}\leq \lambda.$$
Therefore, $\lambda\geq\lambda_{max}=\|X^{T}Y\|_{2}$.

Now we prove the "if" part. If  $\lambda\geq\lambda_{max}$, we can get
$C=\frac{Y}{\lambda}\in \Omega_{D} $.
Under the fact that
$C^{*}(\lambda)=P_{\Omega_{D}}\left(\frac{Y}{\lambda}\right),$
the solution of (3) is
\begin{center}
$C^{*}(\lambda)=\frac{Y}{\lambda}$.
\end{center}
According to the $A^{*}(\lambda)=\lambda C^{*}(\lambda)$ in KKT system (5), we have
$A^{*}(\lambda)=Y$.
By Theorem 2.1, we get
\begin{center}
$\frac{1}{2}\|Y\|^{2}_{F}+\lambda\|B^{*}(\lambda)\|_{*}=\frac{1}{2}\|Y\|^{2}_{F}$.
\end{center}
This yields $\|B^{*}(\lambda)\|_{*}=0$, which implies $B^{*}(\lambda)=0$.
\end{proof}
Before proving Theorem 3.1, we need to  review the basic properties of singular values (Roger 2013).
\begin{Lemma}
Suppose that $P,Q \in R^{p\times q}, l=\rm{min\left\{p, q\right\}}$, two basic inequalities for singular values are
\begin{center}
$\sigma_{i+j-1}\left(P+Q\right)\leq \sigma_{i}\left(P\right)+\sigma_{j}\left(Q\right),  ~~~~~~1 \leq i,j \leq l, i+j \leq l+1;$
\end{center}
\begin{center}
$\sigma_{i+j-1}\left(PQ^{T}\right)\leq \sigma_{i}\left(P\right)\sigma_{j}\left(Q\right), ~~~~~~~~~~~~ 1\leq i,j \leq l, i+j \leq l+1.$
\end{center}
In particular,
\begin{center}
$\sigma_{i}\left(P+Q\right)\leq \sigma_{i}\left(P\right)+\sigma_{1}\left(Q\right),  ~~~~~~1 \leq i \leq l$;
\end{center}
\begin{center}
$\sigma_{1}\left(PQ^{T}\right)\leq \sigma_{1}(P)\sigma_{1}(Q)=\|P\|_{2}\|Q\|_{2}.$
\end{center}
\end{Lemma}
\textbf{Proof of Theorem 3.1}
\begin{proof}
It's known that $C^{*}(\lambda)=P_{\Omega_{D}}(\frac{y}{\lambda})$ from (4). By the nonexpansiveness of $P_{\Omega_{D}}\left(\cdot\right)$, we know that
\begin{center}
$\|C^{*}\left(\lambda\right)-C^{*}\left(\lambda_{0}\right)\|_{F}\leq \left(\frac{1}{\lambda}-\frac{1}{\lambda_{0}}\right)\|Y\|_{F}.$
\end{center}
Setting $\rho:=\left(\frac{1}{\lambda}-\frac{1}{\lambda_{0}}\right)\|Y\|_{F}$, we define the set  $\Omega :=\left\{C\Big|~ \|C-C^{*}\left(\lambda_{0}\right)\|_{F}\leq \rho\right\}$.

In order to prove the desired result, it's enough to consider $\underset{C\in \Omega} \sup ~~ \left\{\sigma_{i}\left(X^{T}C\right)\right\}$ by Theorem 2.2. In fact, if $\underset{C\in \Omega} \sup ~~ \left\{\sigma_{i}\left(X^{T}C\right)\right\}<1$, $\sigma_{i}\left(C^{*}(\lambda)\right)<1$ must holds, which leads to $\sigma_{i}\left(B^{*}(\lambda)\right)=0$.

According to Lemma 5.4, we can get
\begin{align*}
\underset{C\in \Omega} \sup ~~ \left\{\sigma_{i}\left(X^{T}C\right)\right\}
&=\underset{\|D\|_{F}\leq \rho} \sup ~~ \left\{\sigma_{i}\left(X^{T}\left(C^{*}\left(\lambda_{0}\right)+D\right)\right)\right\}\\
&=\underset{\|D\|_{F}\leq \rho} \sup ~~ \left\{\sigma_{i}\left(X^{T}C^{*}\left(\lambda_{0}\right)+X^{T}D\right)\right\}\\
&\leq \underset{\|D\|_{F}\leq \rho} \sup ~~ \left\{\sigma_{i}\left(X^{T}C^{*}\left(\lambda_{0}\right)\right)+\sigma_{1}\left(X^{T}D\right)\right\}\\
&\leq \underset{\|D\|_{F}\leq \rho} \sup ~~ \left\{\sigma_{i}\left(X^{T}C^{*}\left(\lambda_{0}\right)\right)+\|X\|_{2}\|D\|_{2}\right\}\\
&\leq \sigma_{i}\left(X^{T}C^{*}\left(\lambda_{0}\right)\right)+\rho\|X\|_{2}.
\end{align*}
The last inequality is obtained by the fact that $\|D\|_{2}\leq \|D\|_{F}\leq \rho$. Suppose $\sigma_{i}\left(X^{T}C^{*}\left(\lambda_{0}\right)\right)+\rho\|X\|_{2}<1$, that is $\sigma_{i}\left(X^{T}C^{*}\left(\lambda_{0}\right)\right)\leq 1-\left(\frac{1}{\lambda}-\frac{1}{\lambda_{0}}\right)\|X\|_{2}\|Y\|_{F}$. We have
$\underset{C\in \Omega} \sup ~~ \left\{\sigma_{i}\left(X^{T}C\right)\right\}<1$
and $\sigma_{i}\left(X^{T}C^{*}\left(\lambda\right)\right)<1$. Therefore, $\sigma_{i}\left(B^{*}(\lambda)\right)=0$, which implies that for any $j\geq i$, $\sigma_{j}\left(B^{*}(\lambda)\right)=0$. The desired result follows immediately.
\end{proof}
\textbf{Proof of Proposition 3.2}
\begin{proof}
\underline{Case 1}: $X^{T}Y=0$. In this case, $\lambda_{max}=0$ and $B^{*}(\lambda)=0$ for any $\lambda>0$. \\
\underline{Case 2}: All singular values of $X^{T}Y$ are equal to  a nonzero number $\alpha$, which means $\rm rank$$(X^{T}Y)=r$. In this case, $\lambda_{max}=\alpha$ and $C^{*}(\lambda_{max})=\frac{Y}{\alpha}\in \Omega_{D}$. For any $0<\lambda<\alpha$, $$C^{*}(\lambda)=P_{\Omega_{D}}(\frac{Y}{\lambda})=\frac{Y}{\alpha}.$$
Replacing it into the KKT system (5), we know that
$X^{T}Y=X^{T}\left(XB^{*}(\lambda)+\lambda\frac{Y}{\alpha}\right)$, which leads to $$(1-\frac{\lambda}{\alpha})X^{T}Y=X^{T}\left(XB^{*}(\lambda)\right).$$
Hence, $\rm rank$$(X^{T}XB^{*}(\lambda))$=$\rm rank$$(X^{T}Y)=r$. By using the fact that $\rm rank$$(X^{T}XB^{*}(\lambda))\leq$$\rm rank$$(B^{*}(\lambda))\leq r$, we know that $\rm rank$$(B^{*}(\lambda))=r$. Therefore, the result  is proved.
\end{proof}
\textbf{Proof of Theorem 3.2}
\begin{proof}
Under the definition in the theorem, we know that $\lambda_{1}=\lambda_{max}$. The choice of $\lambda$ in Theorem 3.1 should be satisfied that $\lambda<\lambda_{max}$,  so we just talk about the case that $i\geq 2$.  From Theorem 3.1, we can see that when $\lambda>\lambda_{i}$,  rank$\left(B^{*}(\lambda)\right)\leq i-1$. Similarly, if $\lambda>\lambda_{i-1}$, $\rm rank(B^{*}(\lambda))$$\leq i-2$.
Combining these two results, it holds that for any $\lambda \in (\lambda_{i},\lambda_{i-1}]$  $i\geq 2$,
\begin{center}
$\rm rank(B^{*}(\lambda))$$\leq i-1.$
\end{center}
Thus, the conclusion is proved.
\end{proof}
\textbf{Proof of Lemma 3.1}
\begin{proof}
There are two cases in this proof. For every case, the estimate set of $C^{*}\left(\lambda\right)$ can be obtained, then the new interval is contained in $\Omega$ is proved. Before considering these two cases, we give a notation. For any $\lambda \in (0,\lambda_{max}]$ and $t\geq0$, define
\begin{center}
$C_{t}\left(\lambda\right)=C^{*}\left(\lambda\right)+tV_{1}\left(\lambda\right).$
\end{center}
\noindent \underline{Case 1}: $\lambda_{0} \in\left(0,\lambda_{max}\right)$. In this case, $\|X^{T}\left(\frac{Y}{\lambda_{0}}\right)\|_{2}=\frac{1}{\lambda_{0}}\|X^{T}Y\|_{2}
=\frac{\lambda_{max}}{\lambda_{0}}>1$, it causes $\frac{Y}{\lambda_{0}}\notin \Omega_{D}$.

Therefore $\frac{Y}{\lambda_{0}}\neq P_{\Omega_{D}}\left(\frac{Y}{\lambda_{0}}\right)=C^{*}\left(\lambda_{0}\right),$ which leads to
$\frac{Y}{\lambda_{0}}-C^{*}\left(\lambda_{0}\right)\neq0$. By Lemma 5.2 (2), it holds
\begin{eqnarray}
P_{\Omega_{D}}\left(C^{*}\left(\lambda_{0}\right)\right)=C^{*}\left(\lambda_{0}\right).
\end{eqnarray}
A direct result of the definition of $C_{t}(\lambda_{0})$ and Lemma 2.3 is that for any $t\geq 0$,
\begin{eqnarray}
P_{\Omega_{D}}\left(C_{t}\left(\lambda_{0}\right)\right)
=P_{\Omega_{D}}\left(C^{*}\left(\lambda_{0}\right)+t\left(\frac{Y}{\lambda_{0}}-C^{*}\left(\lambda_{0}\right)\right)\right)
=P_{\Omega_{D}}\left(C^{*}\left(\lambda_{0}\right)\right)=C^{*}\left(\lambda_{0}\right).
\end{eqnarray}
Therefore, we have
\begin{align*}
\|C^{*}\left(\lambda\right)-C^{*}\left(\lambda_{0}\right)\|_{F}
&=\left\|P_{\Omega_{D}}\left(\frac{Y}{\lambda}\right)-P_{\Omega_{D}}\left(C_{t}\left(\lambda_{0}\right)\right)\right\|_{F}\\
&\leq \left\|\frac{Y}{\lambda}-C_{t}\left(\lambda_{0}\right)\right\|_{F}\\
&=\left\|t\left(\frac{Y}{\lambda_{0}}-C^{*}(\lambda_{0})\right)-\left(\frac{Y}{\lambda}-C^{*}\left({\lambda_{0}}\right)\right)\right\|_{F}\\
&=\|tV_{1}\left(\lambda_{0}\right)-V_{2}\left(\lambda,\lambda_{0}\right)\|_{F}.
\end{align*}
Because it holds for any $t\geq0$, one can easily see that
\begin{align}
\|C^{*}\left(\lambda\right)-C^{*}\left(\lambda_{0}\right)\|_{F}
&\leq \underset{t} \min \|tV_{1}\left(\lambda_{0}\right)-V_{2}\left(\lambda,\lambda_{0}\right)\|_{F}\\
&=\begin{cases}
\|V_{2}\left(\lambda,\lambda_{0}\right)\|_{F}, & \left\langle V_{1}\left(\lambda_{0}\right),V_{2}\left(\lambda,\lambda_{0}\right)\right\rangle<0\\
\|V_{3}\left(\lambda,\lambda_{0}\right)\|_{F}, & \rm{otherwise}.
\end{cases}
\end{align}
Now we consider the value of  $\left\langle V_{1}\left(\lambda_{0}\right),V_{2}\left(\lambda,\lambda_{0}\right)\right\rangle$. It's easy to see that $0\in \Omega_{D}$. According (7) and the definition of projection operator, we know that
$\|C_{t}\left(\lambda_{0}\right)-C^{*}\left(\lambda_{0}\right)\|^{2}_{F}\leq \|C_{t}\left(\lambda_{0}\right)-0\|^{2}_{F},$
which leads to
$$\|C^{*}\left(\lambda_{0}\right)\|^{2}_{F}-2\left\langle C_{t}\left(\lambda_{0}\right),C^{*}\left(\lambda_{0}\right)\right\rangle\leq 0.$$
Based on the definition of $C_{t}\left(\lambda_{0}\right)$, we have
\begin{eqnarray*}
\|C^{*}\left(\lambda_{0}\right)\|^{2}_{F}
-2\left\langle C^{*}\left(\lambda_{0}\right)+t\left(\frac{Y}{\lambda_{0}}-C^{*}\left(\lambda_{0}\right)\right),C^{*}\left(\lambda_{0}\right)\right\rangle\leq 0.
\end{eqnarray*}
Extending the inner product and transforming all the term to the right side, the last inequality is
\begin{eqnarray*}
0\leq \|C^{*}\left(\lambda_{0}\right)\|^{2}_{F}
+2t\left(\left\langle \frac{Y}{\lambda_{0}},C^{*}\left(\lambda_{0}\right)\right\rangle-\|C^{*}\left(\lambda_{0}\right)\|^{2}_{F}\right).
\end{eqnarray*}
Since the above holds for any $t\geq0$, we obtain
\begin{eqnarray}
\frac{\|Y\|_{F}}{\lambda_{0}}\geq \|C^{*}\left(\lambda_{0}\right)\|_{F}.
\end{eqnarray}
By using the Cauchy-Schwarz inequality, we derive
\begin{align*}
\left\langle V_{1}\left(\lambda_{0}\right),V_{2}\left(\lambda,\lambda_{0}\right)\right\rangle
&=\left\langle\frac{Y}{\lambda_{0}}-C^{*}\left(\lambda_{0}\right),
\frac{Y}{\lambda}-\frac{Y}{\lambda_{0}}
+\frac{Y}{\lambda_{0}}-C^{*}\left(\lambda_{0}\right)\right\rangle\\
&=\left\langle\frac{Y}{\lambda_{0}}-C^{*}\left(\lambda_{0}\right),
\frac{Y}{\lambda}-\frac{Y}{\lambda_{0}}\right\rangle
+\|\frac{Y}{\lambda_{0}}-C^{*}\left(\lambda_{0}\right)\|^{2}_{F}\\
&\geq \left(\frac{1}{\lambda}-\frac{1}{\lambda_{0}}\right)
\left\langle\frac{Y}{\lambda_{0}}-C^{*}\left(\lambda_{0}\right),Y\right\rangle\\
&=\left(\frac{1}{\lambda}-\frac{1}{\lambda_{0}}\right)
\left(\frac{\|Y\|^{2}_{F}}{\lambda_{0}}-\left\langle C^{*}\left(\lambda_{0}\right),Y\right\rangle\right)\\
&\geq \left(\frac{1}{\lambda}-\frac{1}{\lambda_{0}}\right)
\left(\frac{\|Y\|^{2}_{F}}{\lambda_{0}}-\|C^{*}\left(\lambda_{0}\right)\|_{F}\|Y\|_{F}\right).
\end{align*}
This together with (10) yields $\|C^{*}\left(\lambda\right)-C^{*}\left(\lambda_{0}\right)\|_{F}\leq \|V_{3}\left(\lambda,\lambda_{0}\right)\|_{F}$.

Hence, $C^{*}\left(\lambda\right) \in \Omega_{1}$. Replacing $t=1$ into (6), we can get
\begin{align*}
\Omega_{1}&\subseteq \left\{C\Big|~ \|C-C^{*}\left(\lambda_{0}\right)\|_{F}
\leq \|V_{1}\left(\lambda_{0}\right)-V_{2}\left(\lambda,\lambda_{0}\right)\|_{F}\right\}\\
&=\left\{C\Big|~\|C-C^{*}\left(\lambda_{0}\right)\|_{F}
\leq \left\|\frac{Y}{\lambda_{0}}-C^{*}\left(\lambda_{0}\right)-\frac{Y}{\lambda}+C^{*}\left(\lambda_{0}\right)\right\|_{F}\right\}\\
&=\left\{C\Big|~\|C-C^{*}\left(\lambda_{0}\right)\|_{F}
\leq \left(\frac{1}{\lambda}-\frac{1}{\lambda_{0}}\right)\|Y\|_{F}\right\}.
\end{align*}
Therefore, $\Omega_{1}\subseteq \Omega$ in case of $\lambda_{0} \in (0,\lambda_{max})$.\\
\noindent \underline{Case 2}: $\lambda_{0}=\lambda_{max}$. For any $t\geq0$, we want to verify $P_{\Omega_{D}}(C_{t}(\lambda_{max}))=C^{*}(\lambda_{max})$. From Lemma 5.1, we only need to verify that for any $C \in \Omega_{D}$,
\begin{center}
$\left\langle C_{t}\left(\lambda_{max}\right)-C^{*}\left(\lambda_{max}\right),C-C^{*}\left(\lambda_{max}\right)\right\rangle\leq 0$.
\end{center}
\noindent That is, for any $C \in \Omega_{D}$, we need to prove that $\left\langle V_{1}\left(\lambda_{max}\right),C-C^{*}\left(\lambda_{max}\right)\right\rangle\leq 0$. This is true from the fact that
\begin{align*}
&\left\langle V_{1}\left(\lambda_{max}\right),C^{*}\left(\lambda_{max}\right)\right\rangle
=\left\langle\mathcal{V}\left(X\right),\frac{Y}{\lambda_{max}}\right\rangle
=\left\|X^{T}\frac{Y}{\lambda_{max}}\right\|_{2}=1,\\
&\left\langle V_{1}\left(\lambda_{max}\right),C\right\rangle
=\left\langle\mathcal{V}\left(X\right),C\right\rangle
=\|X^{T}C\|_{2}\leq 1.
\end{align*}
Therefore, we have
\begin{align*}
\|C^{*}\left(\lambda\right)-C^{*}\left(\lambda_{max}\right)\|_{F}
&= \left\|P_{\Omega_{D}}\left(\frac{Y}{\lambda}\right)-P_{\Omega_{D}}\left(C_{t}\left(\lambda_{max}\right)\right)\right\|_{F}\\
&\leq \left\|\frac{Y}{\lambda}-C_{t}\left(\lambda_{max}\right)\right\|_{F}\\
&= \left\|tV_{1}\left(\lambda_{max}\right)-\left(\frac{Y}{\lambda}-C^{*}(\lambda_{max})\right)\right\|_{F}\\
&=\left\|tV_{1}\left(\lambda_{max}\right)-V_{2}\left(\lambda,\lambda_{max}\right)\right\|_{F}.
\end{align*}
That is,
$$\|C^{*}\left(\lambda\right)-C^{*}\left(\lambda_{max}\right)\|_{F}
\leq \underset{t\geq0} \min \|tV_{1}\left(\lambda_{max}\right)-V_{2}\left(\lambda,\lambda_{max}\right)\|_{F}
= \|V_{3}\left(\lambda,\lambda_{max}\right)\|_{F}.$$
Hence, $C^{*}\left(\lambda_{max}\right) \in \Omega_{1}$ is proved. Replacing $t=0$ in the last inequality, we can get
\begin{align*}
\Omega_{1}&\subseteq \left\{C\Big|~ \|C-C^{*}\left(\lambda_{max}\right)\|_{F}
\leq \|V_{2}\left(\lambda,\lambda_{max}\right)\|_{F}\right\}\\
&=\left\{C\Big|~\|C-C^{*}\left(\lambda_{max}\right)\|_{F}
\leq \left\|\frac{Y}{\lambda}-C^{*}\left(\lambda_{max}\right)\right\|_{F}\right\}\\
&=\left\{C\Big|~\|C-C^{*}\left(\lambda_{max}\right)\|_{F}
\leq \left(\frac{1}{\lambda}-\frac{1}{\lambda_{max}}\right)\|Y\|_{F}\right\}.
\end{align*}
Thus, $\Omega_{1}\subseteq \Omega$ in the case of $\lambda_{0}=\lambda_{max}$.

From above arguments, we prove $\Omega_{1}\subseteq \Omega$ for any $\lambda_{0}\in (0,\lambda_{max}]$.
\end{proof}
\textbf{Proof of Lemma 3.2}
\begin{proof}
In view of the firm nonexpansiveness of $P_{\Omega_{D}}\left(\cdot\right)$ in Lemma 5.2, we have
\begin{eqnarray}
\left\|C^{*}\left(\lambda\right)-C^{*}\left(\lambda_{0}\right)\right\|^{2}_{F}+
\left\|\frac{Y}{\lambda}-C^{*}\left(\lambda\right)-\frac{Y}{\lambda_{0}}+C^{*}\left(\lambda_{0}\right)\right\|^{2}_{F}
\leq \left\|\frac{Y}{\lambda}-\frac{Y}{\lambda_{0}}\right\|^{2}_{F}.
\end{eqnarray}
This can be reformulated as
\begin{eqnarray*}
\|C^{*}\left(\lambda\right)-C^{*}\left(\lambda_{0}\right)\|^{2}_{F}
\leq \left\langle C^{*}\left(\lambda\right)-C^{*}\left(\lambda_{0}\right),\frac{Y}{\lambda}-\frac{Y}{\lambda_{0}}\right\rangle,
\end{eqnarray*}
which is equivalent to
\begin{eqnarray*}
\left\|C^{*}\left(\lambda\right)-C^{*}\left(\lambda_{0}\right)-\frac{1}{2}\left(\frac{1}{\lambda}-\frac{1}{\lambda_{0}}\right)Y\right\|^{2}_{F}
\leq \frac{1}{4}\left(\frac{1}{\lambda}-\frac{1}{\lambda_{0}}\right)^{2}\|Y\|^{2}_{F}.
\end{eqnarray*}
From the definition of $\Omega_{2}$, we know that $C^{*}\left(\lambda\right)\in\Omega_{2}$ and  $\Omega_{2}\subseteq \Omega$.
\end{proof}
\textbf{Proof of Lemma 3.3}
\begin{proof}
 By using the firm nonexpansiveness of $P_{\Omega_{D}}\left(\cdot\right)$, we have
 \begin{eqnarray*}
 \left\|P_{\Omega_{D}}\left(\frac{Y}{\lambda}\right)-P_{\Omega_{D}}\left(C_{t}\left(\lambda_{0}\right)\right)\right\|^{2}_{F}
 +
\left\|\left(I-P_{\Omega_{D}}\right)\left(\frac{Y}{\lambda}\right)-\left(I-P_{\Omega_{D}}\right)\left(C_{t}\left(\lambda_{0}\right)\right)\right\|^{2}_{F}
\leq
\left\|\frac{Y}{\lambda}-C_{t}\left(\lambda_{0}\right)\right\|^{2}_{F}.
\end{eqnarray*}
where  $C_{t}\left(\lambda\right)$ is defined in the proof of Lemma 3.2. By simple computation and rearranging all terms, the inequality can be transformed into
\begin{eqnarray*}
\left\|P_{\Omega_{D}}\left(\frac{Y}{\lambda}\right)-P_{\Omega_{D}}\left(C_{t}\left(\lambda_{0}\right)\right)\right\|^{2}_{F}
\leq
\left\langle \frac{Y}{\lambda}-C_{t}\left(\lambda_{0}\right), P_{\Omega_{D}}\left(\frac{Y}{\lambda}\right)-P_{\Omega_{D}}\left(C_{t}\left(\lambda_{0}\right)\right)\right\rangle.
\end{eqnarray*}
That is, for any $t\geq0$
\begin{align*}
 \|C^{*}\left(\lambda\right)-C^{*}\left(\lambda_{0}\right)\|^{2}_{F}
 & \leq \left\langle \frac{Y}{\lambda}-C_{t}\left(\lambda_{0}\right), C^{*}\left(\lambda\right)-C^{*}\left(\lambda_{0}\right)\right\rangle\\
& =\left\langle \frac{Y}{\lambda}-C^{*}\left(\lambda_{0}\right)-tV_{1}\left(\lambda_{0}\right), C^{*}\left(\lambda\right)-C^{*}\left(\lambda_{0}\right)\right\rangle\\
&=\left\langle V_{2}\left(\lambda,\lambda_{0}\right)-tV_{1}\left(\lambda_{0}\right), C^{*}\left(\lambda\right)-C^{*}\left(\lambda_{0}\right)\right\rangle.
 \end{align*}
By transforming all terms to the left side, the equality can be written as
 \begin{align*}
 \|C^{*}\left(\lambda\right)-C^{*}\left(\lambda_{0}\right)\|^{2}_{F}-\left\langle V_{2}\left(\lambda,\lambda_{0}\right)-tV_{1}\left(\lambda_{0}\right), C^{*}\left(\lambda\right)-C^{*}\left(\lambda_{0}\right)\right\rangle=0,
 \end{align*}
which leads to
\begin{eqnarray}
\left\|C^{*}\left(\lambda\right)-C^{*}\left(\lambda_{0}\right)
 -\frac{1}{2}\left(V_{2}\left(\lambda,\lambda_{0}\right)-tV_{1}\left(\lambda_{0}\right)\right)\right\|^{2}_{F}
 \leq \frac{1}{4}\|V_{2}\left(\lambda,\lambda_{0}\right)-tV_{1}\left(\lambda_{0}\right)\|^{2}_{F}.
 \end{eqnarray}
According to the proof of Lemma 3.1, we have
$$\frac{\left\langle V_{1}\left(\lambda_{0}\right),V_{2}\left(\lambda,\lambda_{0}\right)\right\rangle}
 {\|V_{1}\left(\lambda_{0}\right)\|^{2}_{F}}\geq 0.$$
Define $t$ as above and replace it into (12), we have
\begin{center}
$\left\|C^{*}\left(\lambda\right)-C^{*}\left(\lambda_{0}\right)
-\frac{1}{2}V_{3}\left(\lambda,\lambda_{0}\right)\right\|_{F}
\leq \frac{1}{2} \|V_{3}\left(\lambda,\lambda_{0}\right)\|_{F}$.
\end{center}
\end{proof}

\end{document}